\documentclass[a4paper]{article}
\usepackage{graphicx}
\usepackage{amsmath}
\usepackage{amsfonts}
\usepackage{amssymb}
\newtheorem{theorem}{Theorem}[subsection]

\newtheorem{conjecture}[theorem]{Conjecture}
\newtheorem{corollary}[theorem]{Corollary}

\newtheorem{definition}[theorem]{Definition}
\newtheorem{example}[theorem]{Example}

\newtheorem{lemma}[theorem]{Lemma}
\newtheorem{notation}[theorem]{Notation}

\newtheorem{proposition}[theorem]{Proposition}
\newtheorem{remark}[theorem]{Remark}

\newenvironment{proof}[1][Proof]{\textbf{#1.} }{\ \rule{0.5em}{0.5em}}

\begin{document}

\title{Schensted-Type correspondence, Plactic Monoid and Jeu de Taquin for type
$C_{n}$}
\author{C\'{e}dric Lecouvey\\lecouvey@math.unicaen.fr}
\date{}
\maketitle
\begin{abstract}
We use Kashiwara's theory of crystal bases to study the plactic monoid for
$U_{q}(sp_{2n})$. Then we describe the corresponding insertion and sliding
algorithms. The sliding algorithm is essentially the symplectic Jeu de Taquin
defined by Sheats and our construction gives the proof of its compatibility
with plactic relations.
\end{abstract}

\section{Introduction}

It is well known that the Schensted bumping algorithm yields a bijection
between words $w$ of length $l$ on the ordered alphabet $\mathcal{A}%
_{n}=\{1<2<\cdot\cdot\cdot<n\}$ and pairs $(P(w),Q(w))$ of tableaux of the
same shape containing $l$ boxes where $P(w)$ is a \ semi-standard Young
tableau on $\mathcal{A}_{n}$ and $Q(w)$ is a standard tableau. This bijection
is called the Robinson-Schensted correspondence (see e.g.\cite{Fu}). Notice
that the tableau $P(w)$ may also be constructed from $w$ by using the
Sch\"{u}tzenberger sliding algorithm.

We can define a relation $\sim$ on the free monoid $\mathcal{A}_{n}^{\ast}$
by:
\[
w_{1}\sim w_{2}\Longleftrightarrow P(w_{1})=P(w_{2}).
\]
Then the quotient $Pl(A_{n}):=\mathcal{A}_{n}^{\ast}/\sim$ can be described as
the quotient of $\mathcal{A}_{n}^{\ast}$ by Knuth relations:
\begin{align*}
zxy  &  =xzy\text{ \ \ and \ \ }yzx=yxz\text{ if }x<y<z,\\
xyx  &  =xxy\text{ \ \ and \ \ }xyy=yxy\text{ if }x<y.
\end{align*}
Hence $Pl(A_{n}),$ which may be identified with the set of standard Young
tableaux, becomes a monoid in a natural way. This monoid is called the
''plactic monoid'' and has been introduced by Lascoux and Sch\"{u}tzenberger
in order to give an illuminating proof of the Littlewood-Richardson rule for
decomposing tensor products of irreducible $gl_{n}$-modules \cite{LS}.

There have been attempts to find a Robinson-Schensted type correspondence and
plactic relations for the other classical Lie algebras. In \cite{Be}, Berele
has explained an insertion algorithm for $sp_{2n}$ and in \cite{SU} Sundaram
gives an insertion scheme for $so_{2n+1}$ but it seems difficult to obtain
plactic relations from these schemes. More recently Littelmann has used his
path model to introduce a plactic algebra for any simple Lie algebra
\cite{Lit}. In the symplectic case, a set of defining relations of the plactic
algebra was described independently (without proof) by Lascoux, Leclerc and
Thibon \cite{LLT}. In \cite{SH}, Sheats has introduced a symplectic Jeu de
Taquin analogous to Sch\"{u}tzenberger's sliding algorithm and has conjectured
its compatibility with the plactic relations of \cite{LLT}.

The Robinson-Schensted correspondence has a natural interpretation in terms of
Kashiwara's theory of crystal bases \cite{DJM}, \cite{Ka2}, \cite{LLT}. Let
$W_{n}$ denote the vector representation of $gl_{n}$. By considering each
vertex of the crystal graph of $\underset{l\geq0}{\bigoplus}W_{n}^{\bigotimes
l}$ as a word on $\mathcal{A}_{n}$, we have for any words $w_{1}$ and $w_{2}$:

\begin{itemize}
\item $P(w_{1})=P(w_{2})$ if and only if $w_{1}$ and $w_{2}$ occur at the same
place in two isomorphic connected components of this graph.

\item $Q(w_{1})=Q(w_{2})$ if and only if $w_{1}$ and $w_{2}$ occur in the same
connected component of this graph.
\end{itemize}

Replacing $W_{n}$ by the vector representation $V_{n}$ of $sp_{2n}$ whose
basis vectors are labelled by
\[
\mathcal{C}_{n}=\{1<\cdot\cdot\cdot<n-1<n<\overline{n}<\overline{n-1}%
<\cdot\cdot\cdot<\overline{1}\},
\]
we can define similarly a relation on $\mathcal{C}_{n}^{\ast}$ by: $w_{1}%
\sim_{n}w_{2}$ if and only if $w_{1}$ and $w_{2}$ have the same position in
two isomorphic connected components of the crystal of $G_{n}=\underset
{l}{\bigoplus}V_{n}^{\bigotimes l}$ and set $Pl(C_{n}):=\mathcal{C}_{n}^{\ast
}/\sim_{n}$. In this article, we undertake a detailed investigation of
$Pl(C_{n})$ and of the corresponding insertion and sliding algorithms.

We first recall in part 2 the combinatorial notion of symplectic tableau
introduced by De Concini \cite{DeC} (analogous to Young tableaux for type
$C_{n}$) and how it relates to Kashiwara's theory of crystal graphs. In part 3
we derive a set of defining relations for $Pl(C_{n})$ similar to Knuth
relations using the description of the crystal graphs given by Kashiwara and
Nakashima \cite{KN}. Our plactic relations are those of \cite{LLT} up to a
small mistake that we correct. We will recover Littelmann's interpretation
\cite{Lit} of plactic relations in terms of crystal isomorphisms but not his
explicit description of $Pl(C_{n})$ which is not totally exact.\ In part 4, we
describe a column insertion algorithm analogous to the bumping algorithm for
type $A$. We note that the insertion procedure is more complicated as for type
A and cannot be described as a mere bumping algorithm. Using the notion of
oscillating tableaux (analogous to standard tableaux for type $C_{n}$), this
algorithm yields the desired Robinson-Schensted type correspondence in part 5.
Finally we recall in part 6 the definition of the symplectic Jeu de Taquin
introduced by Sheats. Then we prove two conjectures of Sheats: applying
$\mathrm{SJDT}$ slides stays within the initial plactic congruence class and
the result of the rectification of any skew symplectic tableau is independent
of the order in which the inner corners are filled. In particular Sheats Jeu
de Taquin can be used to compute the $P$ symbol of a word $w\in\mathcal{C}%
_{n}^{\ast}$.

Similar results for types $B$ and $D$ will be discussed in a forthcoming paper
\cite{Lec}.

\bigskip

\noindent\textbf{Note }\textit{While writing this work, I have been informed
that T. H. Baker \cite{ba} has obtained independently and by different methods
essentially the same insertion schemes as those described in Section 4.}

\section{Crystal graphs and symplectic tableaux}

This section focuses on the notion of symplectic tableaux introduced by
Kashiwara and Nakashima to label the vertices of the crystal graphs of a
$U_{q}(sp_{2n})$-module \cite{KN}. We relate the symplectic tableaux of
Kashiwara and Nakashima to the tableaux used by De Concini in \cite{DeC} to
express the irreducible characters of the classical group $Sp_{2n}$.

\subsection{Convention for crystal graphs of $U_{q}(sp_{2n})$-modules}

We adopt Kashiwara's convention \cite{Ka2} for crystal graphs. The Dynkin
diagram of $sp_{2n}$ is labelled by:
\[
\overset{1}{\circ}-\overset{2}{\circ}-\overset{3}{\circ}\cdot\cdot
\cdot\overset{n-2}{\circ}-\overset{n-1}{\circ}\Longleftarrow\overset{n}{\circ}%
\]
Accordingly, the crystal graph of the vector representation $V_{n}$ of
$U_{q}(sp_{2n})$ will be labelled as follows:
\[
1\overset{1}{\rightarrow}2\cdot\cdot\cdot\cdot\rightarrow n-1\overset
{n-1}{\rightarrow}n\overset{n}{\rightarrow}\overline{n}\overset{n-1}%
{\rightarrow}\overline{n-1}\overset{n-2}{\rightarrow}\cdot\cdot\cdot
\cdot\rightarrow\overline{2}\overset{1}{\rightarrow}\overline{1}%
\]
Recall that crystal graphs are oriented colored graphs with colors
$i\in\{1,...,n\}.$ An arrow $a\overset{i}{\rightarrow}b$ means that
$\widetilde{f}_{i}(a)=b$ and $\widetilde{e}_{i}(b)=a$ where $\widetilde{e}%
_{i}$ and $\widetilde{f}_{i}$ are the crystal graph operators (for a review of
crystal bases and crystal graphs see \cite{Ka2}). The action of $\widetilde
{e}_{i}$ and $\widetilde{f}_{i}$ on the tensor product of two crystal graphs
(which is the crystal graph of the tensor product of these representations) is
given by :%

\begin{align}
\widetilde{f_{i}}(u%
{\textstyle\bigotimes}
v)  &  =\left\{
\begin{tabular}
[c]{c}%
$\widetilde{f}_{i}(u)\bigotimes v$ if $\varphi_{i}(u)>\varepsilon_{i}(v)$\\
$u\bigotimes\widetilde{f}_{i}(v)$ if $\varphi_{i}(u)\leq\varepsilon_{i}(v)$%
\end{tabular}
\right. \label{tens}\\
\widetilde{e_{i}}(u%
{\textstyle\bigotimes}
v)  &  =\left\{
\begin{tabular}
[c]{c}%
$u\bigotimes\widetilde{e_{i}}(v)$ if $\varphi_{i}(u)<\varepsilon_{i}(v)$\\
$\widetilde{e_{i}}(u)\bigotimes v$ if $\varphi_{i}(u)\geq\varepsilon_{i}(v)$%
\end{tabular}
\right. \nonumber
\end{align}
where $\varepsilon_{i}(u)=\max\{k;\widetilde{e}_{i}^{k}(u)\neq0\}$ and
$\varphi_{i}(u)=\max\{k;\widetilde{f}_{i}^{k}(u)\neq0\}$. By induction, this
allows us to define a crystal graph for the representations $V_{n}^{\bigotimes
l}$ for any $l$. Each vertex $u_{1}\bigotimes u_{2}\bigotimes\cdot\cdot
\cdot\bigotimes u_{l}$ of the crystal graph of $V_{n}^{\bigotimes l}$ will be
identified with the word $w=u_{1}u_{2}\cdot\cdot\cdot u_{l}$ over the finite
totally ordered alphabet
\[
\mathcal{C}_{n}=\{1<\cdot\cdot\cdot<n<\overline{n}<\cdot\cdot\cdot
<\overline{1}\}.
\]
Write $\mathcal{C}_{n}^{\ast}$ for the free monoid defined on $\mathcal{C}%
_{n}$. A barred (resp. unbarred) letter is a letter of $\mathcal{C}_{n}$ $>n$
(resp. $<\overline{n}$). We denote by $d(w)=(d_{1},...,d_{n})$ the $n$-tuple
where $d_{i}$ is the number of letters $i$ in $w$ minus the number of letters
$\overline{i}$, and by $\frak{L}(w)$ the length of this word.

Let $G_{n}$ and $G_{n,l}$ be the crystal graphs of $T(V_{n})=\underset
{l}{\bigoplus}V_{n}^{\bigotimes l}$ and $V_{n}^{\bigotimes l}$. Then the
vertices of $G_{n}$ are indexed by the words of $\mathcal{C}_{n}^{\ast}$ and
those of $G_{n,l}$ by the words of length $l$. $G_{n,l}$ may be decomposed
into connected components which are the crystal graphs of the irreducible
representations occurring in the decomposition of the representation
$V_{n}^{\bigotimes l}$. If $w$ is vertex of $G_{n}$, we denote by $B(w)$ the
connected component of $G_{n}$ containing $w$. In the sequel we call
sub-crystal of $G_{n}$ a union of connected components of $G_{n}$.\ Two
sub-crystals $\Gamma$ and $\Gamma^{\prime}$ of $G_{n}$ are isomorphic when
there exists a bijective map $\xi$ from $\Gamma\cup\{0\}$ to $\Gamma^{\prime
}\cup\{0\}$ such that $\xi(0)=0$ and%
\[
\xi\circ\widetilde{e}_{i}=\widetilde{e}_{i}\circ\xi\text{, }\xi\circ
\widetilde{f}_{i}=\widetilde{f}_{i}\circ\xi\text{ for }i=1,...,n.
\]
Then we will say that $\xi$ is a crystal isomorphism.\ We now introduce the
coplactic relation.

\begin{definition}
Let $w_{1}$ and $w_{2}\in$ $\mathcal{C}_{n}^{\ast}$. We write $w_{1}%
\longleftrightarrow w_{2}\ $if and only if $w_{1}$ and $w_{2}$ belong to the
same connected component of $G_{n}.$
\end{definition}

\noindent Note that $w_{1}\longleftrightarrow w_{2}$ if and only if
$w_{2}=\widetilde{H}(w_{1})$ where $\widetilde{H}$ is a product of Kashiwara's
operators for $U_{q}(sp_{2n})$. The Lemma below follows immediately from
(\ref{tens}).

\begin{lemma}
\label{lem_coplactic}If $w_{1}=u_{1}v_{1}$ and $w_{2}=u_{2}v_{2}$ with
$\frak{L}(u_{1})=\frak{L}(u_{2})$ and $\frak{L}(v_{1})=\frak{L}(v_{2})$%
\[
w_{1}\longleftrightarrow w_{2}\Longrightarrow\left\{
\begin{tabular}
[c]{l}%
$u_{1}\longleftrightarrow u_{2}$\\
$v_{1}\longleftrightarrow v_{2}$%
\end{tabular}
\right.  .
\]
\end{lemma}

Each connected component contains a unique vertex $v^{0}$ such that
$\widetilde{e}_{i}(v^{0})=0$ for $i=1,...,n$. We will call it the vertex of
highest weight. Let $\Lambda_{1},...,\Lambda_{n}$ be the fundamental weights
of $U_{q}(sp_{2n})$. The highest weight of the vector representation $V_{n}$
is equal to $\Lambda_{1}$. If $v^{0}$ is the vertex of highest weight of the
crystal graph of an irreducible $U_{q}(sp_{2n})$-module contained in
$V_{n}^{\bigotimes l}$, the weight $wt(v^{0})$ of $v^{0}$ is the highest
weight of this representation. It is given by:
\[
\mathrm{wt}(v^{0})=d_{n}\Lambda_{n}+\overset{n-1}{\underset{i=1}{\sum}}%
(d_{i}-d_{i+1})\Lambda_{i}%
\]
where $d(v^{0})=(d_{1},...,d_{n})$. Then two connected components $B$ and
$B^{\prime}$ are isomorphic if and only if their vertices of highest weight
have the same weight. Note that in this case the isomorphism between $B$ and
$B^{\prime}$ is unique.\ This implies the following lemma that we will need in
Section \ref{Sec_JDT}.

\begin{lemma}
\label{lem_iso_G}Consider $\Gamma$ be a sub-crystal of $G_{n}$. Denote by
$\Gamma^{0}$ the set of highest weight vertices of $\Gamma$. Let $\xi$ be a
map $\Gamma\cup\{0\}\rightarrow G_{n}$ such that

\noindent$\mathrm{(i)}$ $\xi(0)=0,$

\noindent$\mathrm{(ii)}$ if $w^{0}\in\Gamma^{0},$ $\xi(w^{0})$ is a highest
weight vertex, $d(\xi(w^{0}))=d(w^{0})$ and the restriction of $\xi$ to
$\Gamma^{0}$ is injective$,$

\noindent$\mathrm{(iii)}$ for $w\in\Gamma$, and $i=1,...,n$ such that
$\widetilde{f}_{i}(w)\neq0$, $\xi(\widetilde{f}_{i}w)=\widetilde{f}_{i}\xi(w)$.

Then $\xi$ is a crystal isomorphism from $\Gamma$ to $\xi(\Gamma)$.
\end{lemma}

It is convenient to parametrize the irreducible $U_{q}(sp_{2n})$-modules by
Young diagrams. Let $\lambda=\underset{i=1}{\overset{n}{\sum}}\lambda
_{i}\Lambda_{i}$ with $\lambda_{i}\in\mathbb{N}$. We associate to $\lambda$
the Young diagram $Y_{\lambda}$ containing exactly $\lambda_{i}$ columns of
height $i$ and we set $\left|  \lambda\right|  =\overset{n-1}{\underset
{i=0}{\sum}}\lambda_{i}\,i$ (i.e. $\left|  \lambda\right|  $ is the number of
boxes of $Y_{\lambda}$). In the sequel we say that $Y_{\lambda}$ has shape
$\lambda$. We write $B(\lambda)$ for the crystal graph of $V(\lambda)$, the
irreducible module of highest weight $\lambda$.

The following lemma comes from formula (\ref{tens}):

\begin{lemma}
\label{lmmaximalweight}For any words $w_{1}$ and $w_{2}$ in $\mathcal{C}%
_{n}^{\ast}$, the word $w_{1}w_{2}$ is a vertex of highest weight of a
connected component of $G_{n}$ if and only if:

\begin{itemize}
\item $w_{1}$ is a vertex of highest weight $($i.e. $\varepsilon_{i}(w_{1})=0$
for $i=1,...,n)$

\item  for any $i=1,...,n$ we have $\varepsilon_{i}(w_{2})\leq\varphi
_{i}(w_{1}).$
\end{itemize}
\end{lemma}

We now recall a simple process described by Kashiwara and Nakashima \cite{KN}
to compute the action of the crystal graphs operators $\widetilde{e}_{i}$ and
$\widetilde{f}_{i}$ on a word $w$ (this is a simple consequence of
(\ref{tens})). The idea is to consider first the subword $w_{i}$ of $w$
containing only the letters $\overline{i+1},\overline{i},i,i+1$. Then we
encode in $w_{i}$ each letter $\overline{i+1}$ or $i$ by the symbol $+$ and
each letter $\overline{i}$ or $i+1$ by the symbol $-$. Since $\widetilde
{e}_{i}(+-)=\widetilde{f}_{i}(+-)=0$ the factors of type $+-$ may be ignored
in $w_{i}.$ So we obtain a subword $w_{i}^{(1)}$ in which we can ignore all
the factors $+-$ to construct a new subword $w_{i}^{(2)}$ etc... Finally we
obtain a subword
\begin{equation}
\rho_{i}(w)=-^{r}+^{s}. \label{regle+-}%
\end{equation}

\begin{itemize}
\item  if $r>0,\widetilde{e}_{i}(w)$ is obtained by changing the rightmost
symbol $-$ of $\rho(w)$ into its corresponding symbol $+$ (i.e. $i+1$ into $i$
and $\overline{i}$ into $\overline{i+1}$) the others letters of $w$ being
unchanged. If $r=0,\widetilde{e}_{i}(w)=0.$

\item  if $s>0,\widetilde{f}_{i}(w)$ is obtained by changing the leftmost
symbol $+$\ of $\rho(w)$ into its corresponding symbol $-$ (i.e. $i$ into
$i+1$ and $\overline{i+1}$ into $\overline{i}$) the others letters of $w$
being unchanged. If $s=0,\widetilde{f}_{i}(w)=0.$
\end{itemize}

\begin{example}
If
\[
w=12\overline{1}3\overline{1}\,\overline{2}21211\overline{1}\,\overline
{3}21\overline{1}\,\overline{2}\,\overline{2}42
\]
and $i=1$ we have%
\[
w_{i}=\mathbf{12}\overline{1}\,\overline{1}\,\overline{\mathbf{2}}%
\mathbf{212}1\mathbf{1}\overline{\mathbf{1}}2\mathbf{1}\overline{\mathbf{1}%
}\,\overline{2}\,\overline{\mathbf{2}}\mathbf{2},\text{ \ \ \ }w_{i}%
^{(1)}=\overline{1}\,\overline{1}\,\mathbf{12}\,\overline{2}\,,\text{
\ \ \ }w_{i}^{(2)}=\overline{1}\,\overline{1}\,\overline{2}%
\]%
\[
\rho_{i}(w)=--+
\]%
\[
\widetilde{e}_{i}(w)=12\overline{1}3\overline{2}\,\overline{2}21211\overline
{1}\,\overline{3}21\overline{1}\,\overline{2}\,\overline{2}42,\text{
}\widetilde{f}_{i}(w)=12\overline{1}3\overline{1}\,\overline{2}21211\overline
{1}\,\overline{3}21\overline{1}\,\overline{1}\,\overline{2}42.
\]
\end{example}

\subsection{Admissible and coadmissible columns \label{subsec_splitC}}

A column is a Young diagram $C$ of column shape filled with letters of
$\mathcal{C}_{n}$ strictly increasing from top to bottom. The word of
$\mathcal{C}_{n,}^{\ast}$obtained by reading the letters of $C$ from top to
bottom is called the reading of $C$ and denoted by \textrm{w(}$C),$ that is,
we write%
\[
\text{\textrm{w}}(C)=x_{1}\cdot\cdot\cdot x_{l}\text{ for }C=%
\begin{tabular}
[c]{|l|}\hline
$x_{1}$\\\hline
$\cdot$\\\hline
$\cdot$\\\hline
$\cdot$\\\hline
$x_{l}$\\\hline
\end{tabular}
\text{ with }x_{1}<\cdot\cdot\cdot<x_{l}\text{.}%
\]
The height of $C$ denoted by $h(C)$ is the number of letters in $C$. A word
$w$ is a column word if there exists a column $C$ such that $w=$%
\textrm{w(}$C)$, i.e. if $w$ is strictly increasing.

\begin{definition}
Let $C$ be a column such that \textrm{w(}$C)=x_{1}\cdot\cdot\cdot x_{h(C)}%
$.\ Then $C$ is KN-admissible if there is no pair $(z,\bar{z})$ of letters in
$C$ such that:
\[
z=x_{p},\overline{z}=x_{q}\text{ and }q-p<h(C)-z+1
\]
\end{definition}

\begin{remark}
\label{remNz}The maximal height of a KN-admissible column is $n$. Moreover,
$C$ is KN-admissible if and only if, for any $m=1,...,n$ the number $N(m)$ of
letters $x$ in $C$ such that $x\leq m$ or $x\geq\overline{m}$ satisfies
$N(m)\leq m.$ Moreover if there exists in $C$ a letter $m\leq n$ such that
$N(m)>m$ then $C$ contains a pair $(z,\overline{z})$ satisfying $N(z)>z$.
\end{remark}

We will need a different version of admissible columns, which goes back to
\cite{DeC}. The equivalence between the two notions is proved in \cite{SH}.

\begin{definition}
Let $C$ be a column and $I=\{z_{1}>\cdot\cdot\cdot>z_{r}\}$ the set of
unbarred letters $z$ such that the pair $(z,\overline{z})$ occurs in $C$. The
column $C$ can be split when there exists (see the example below) a set of $r$
unbarred letters $J=\{t_{1}>\cdot\cdot\cdot>t_{r}\}\subset\mathcal{C}_{n}$
such that:

$\ \ \ t_{1}$ is the greatest letter of $\mathcal{C}_{n}$ satisfying:
$t_{1}<z_{1},t_{1}\notin C$ and $\overline{t_{1}}\notin C,$

\ \ \ for $i=2,...,r$, $t_{i}$ is the greatest letter of $\mathcal{C}_{n}$
satisfying: $t_{i}<\min(t_{i-1,}z_{i}),$ $t_{i}\notin C$ and $\overline{t_{i}%
}\notin C.$

\noindent In this case we write:

\noindent$rC$ for the column obtained by changing in $C,$ $\overline{z}_{i}$
into $\overline{t}_{i}$ for each letter $z_{i}\in I$ and by reordering if necessary,

\noindent$lC$ for the column obtained by changing in $C,$ $z_{i}$ into $t_{i}$
for each letter $z_{i}\in I$ and by reordering if necessary.
\end{definition}

\begin{proposition}
(Sheats) A column $C$ is KN-admissible if and only if it can be split.
\end{proposition}

\begin{example}
Let $C=2467\overline{7}\,\overline{4}\,\overline{2}$. Then:
\[
I=\{7,4,2\};\text{ }J=\{5,3,1\};\text{ }lC=1356\overline{7}\,\overline
{4}\,\overline{2}\text{ and }rC=2467\overline{5}\,\overline{3}\,\overline
{1}\,.
\]
Therefore $C$ can be split. Notice that $C^{\prime}=24567\overline
{7}\,\overline{4}\,\overline{2}$ can not be split. Indeed we have
$I_{C^{\prime}}=\{7,4,2\}$. Then $t_{1}=3,$ $t_{2}=1$ (with the notation of
the above definition) but it is impossible to find $t_{3}<1$ in $\mathcal{C}%
_{n}$.
\end{example}

An admissible column word is the reading of an admissible column.\ In Section
5 we will need the notion of KN-coadmissible column. A column $C$ is said to
be KN-coadmissible if for each pair $(z,\overline{z})$ in $C$, the number
$N^{\ast}(z)$ of letters $x$ in $C$ such that $x\geq z$ or $x\leq\overline{z}$
satisfies
\begin{equation}
N^{\ast}(z)\leq n-z+1. \label{N_star}%
\end{equation}
Let $C$ be a KN-admissible column. Denote by $C^{\ast}$ the column obtained by
filling the shape of $C$ (from top to bottom) with the unbarred letters of
$lC$ in increasing order followed by the barred letters of $rC$ in increasing
order. Then it is easy to prove that $C^{\ast}$ is KN-coadmissible. More
precisely the map:%
\begin{equation}
\Phi:C\mapsto C^{\ast} \label{phi}%
\end{equation}
is a bijection between the sets of KN-admissible and KN-coadmissible columns
of the same height. Starting from a KN-coadmissible column $C^{\ast}$ we can
compute the pair $(lC,rC)$ associated to the unique KN-admissible column $C$
such that $C^{\ast}=C$ by reversing the previous algorithm. Then $C$ is the
column containing the unbarred letters of $rC$ and the barred letters of $lC$.

\begin{example}
\ \ \ \ \ \ \ \ \ \ 

if $C=;%
\begin{tabular}
[c]{|l|}\hline
$\mathtt{1}$\\\hline
$\mathtt{4}$\\\hline
$\mathtt{\bar{5}}$\\\hline
$\mathtt{\bar{4}}$\\\hline
$\mathtt{\bar{3}}$\\\hline
\end{tabular}
$ then $(lC,rC)=\left(  \text{ }%
\begin{tabular}
[c]{|l|}\hline
$\mathtt{1}$\\\hline
$\mathtt{2}$\\\hline
$\mathtt{\bar{5}}$\\\hline
$\mathtt{\bar{4}}$\\\hline
$\mathtt{\bar{3}}$\\\hline
\end{tabular}
,%
\begin{tabular}
[c]{|l|}\hline
$\mathtt{1}$\\\hline
$\mathtt{4}$\\\hline
$\mathtt{\bar{5}}$\\\hline
$\mathtt{\bar{3}}$\\\hline
$\mathtt{\bar{2}}$\\\hline
\end{tabular}
\text{ }\right)  $ and $C^{\ast}=%
\begin{tabular}
[c]{|l|}\hline
$\mathtt{1}$\\\hline
$\mathtt{2}$\\\hline
$\mathtt{\bar{5}}$\\\hline
$\mathtt{\bar{3}}$\\\hline
$\mathtt{\bar{2}}$\\\hline
\end{tabular}
.\vspace{0.2cm}$

if $D^{\ast}=%
\begin{tabular}
[c]{|l|}\hline
$\mathtt{1}$\\\hline
$\mathtt{2}$\\\hline
$\mathtt{\bar{2}}$\\\hline
$\mathtt{\bar{1}}$\\\hline
\end{tabular}
;$ then $(lD,rD)=\left(  \text{ }%
\begin{tabular}
[c]{|l|}\hline
$\mathtt{1}$\\\hline
$\mathtt{2}$\\\hline
$\mathtt{\bar{4}}$\\\hline
$\mathtt{\bar{3}}$\\\hline
\end{tabular}
,%
\begin{tabular}
[c]{|l|}\hline
$\mathtt{3}$\\\hline
$\mathtt{4}$\\\hline
$\mathtt{\bar{2}}$\\\hline
$\mathtt{\bar{1}}$\\\hline
\end{tabular}
\text{ }\right)  $ and $D=%
\begin{tabular}
[c]{|l|}\hline
$\mathtt{3}$\\\hline
$\mathtt{4}$\\\hline
$\mathtt{\bar{4}}$\\\hline
$\mathtt{\bar{3}}$\\\hline
\end{tabular}
.$
\end{example}

The original definition of admissible and coadmissible columns due to De
Concini and used in \cite{DeC} and \cite{SH} differs from the above one in
that De Concini uses the alphabet $\{\overline{n}<\overline{n-1}<\cdot
\cdot\cdot<\overline{1}<1<\cdot\cdot\cdot<n-1<n\}$ instead of $\mathcal{C}%
_{n}$. Up to this change of notation, De Concini coadmissible columns provide
a natural labelling of the crystal graphs $B(\Lambda_{p})$ $p=1,...,n$. The
one-to-one correspondence between De Concini's admissible columns and
KN-admissible columns is explicitly described by $\Phi$ composed with the
permutation $p\longleftrightarrow\overline{n-p+1}$, $\overline{p}%
\longleftrightarrow n-p+1$. In the sequel we systematically translate from De
Concini's convention to Kashiwara's convention. To simplify the notation we
will write admissible and coadmissible for KN-admissible and KN-coadmissible.

\subsection{Symplectic tableaux}

\begin{definition}
(Kashiwara-Nakashima \cite{KN}).\ Let $C_{1}$ and $C_{2}$ be admissible
columns with $h(C_{2})\leq h(C_{2})$ and consider the tableau%
\[
C_{1}C_{2}=%
\begin{tabular}
[c]{|l|l}\hline
$i_{1}$ & \multicolumn{1}{|l|}{$j_{1}$}\\\hline
$\cdot$ & \multicolumn{1}{|l|}{$\cdot$}\\\hline
$\cdot$ & \multicolumn{1}{|l|}{$\cdot$}\\\hline
$\cdot$ & \multicolumn{1}{|l|}{$j_{M}$}\\\hline
$i_{N}$ & \\\cline{1-1}%
\end{tabular}
.
\]
For $1\leq a\leq b\leq n$, $C_{1}C_{2}$ contains an $(a,b)$-configuration if
there exists $1\leq p\leq q<r\leq s\leq M$ such that $i_{p}=a,j_{q}%
=b,j_{r}=\overline{b},j_{s}=\overline{a}$ or $i_{p}=a,i_{q}=b,i_{r}%
=\overline{b},j_{s}=\overline{a}$. Then we denote by $p(a,b)$ the positive
integer:%
\[
p(a,b)=(s-r)+(q-p).
\]
\end{definition}

\begin{definition}
(Kashiwara-Nakashima \cite{KN}). Let $C_{1},C_{2},...,C_{r}$ be admissible
columns such that $h(C_{i})\geq h(C_{i+1})$ for $i=1,...,r-1$. $T=C_{1}%
C_{2}\cdot\cdot\cdot C_{r}$ is a symplectic tableau when its rows are weakly
increasing from left to right and for $i=1,...,r-1$, $C_{i}C_{i+1}$ does not
contain an $(a,b)$-configuration with $p(a,b)\geq b-a$.
\end{definition}

Let $C_{1}$ and $C_{2}$ be two admissible columns. We write:

\begin{itemize}
\item $C_{1}\leq C_{2}$ when $h(C_{1})\geq h(C_{2})$ and the rows of the
tableau $C_{1}C_{2}$ are weakly increasing from left to right.

\item $C_{1}\preceq C_{2}$ when $rC_{1}\leq lC_{2}$.
\end{itemize}

Note that for any admissible column $C$, $lC\leq C\leq rC$, hence
$C_{1}\preceq C_{2}\Longrightarrow C_{1}\leq C_{2}$.\ The next proposition is
proved in \cite{SH} (Lemma A.4) with the convention of De Concini. It gives a
description of the symplectic tableaux which is the translation into
Kashiwara's convention of the description used by Sheats in \cite{SH}.

\begin{proposition}
Let $C_{1},C_{2},...,C_{r}$ be admissible columns. $T=C_{1}C_{2}\cdot
\cdot\cdot C_{r}$ is a symplectic tableau if and only if $C_{i}\preceq
C_{i+1}$ for $i=1,...,r-1.$
\end{proposition}

Let $T=C_{1}C_{2}\cdot\cdot\cdot C_{r}$ be a symplectic tableau. The shape of
$T$ is the shape of the Young diagram associated to $T$.\ The reading of $T$
is the word \textrm{w(}$T)$ defined by $\mathrm{w(}T)=\mathrm{w(}%
C_{r})\mathrm{w(}C_{r-1})\cdot\cdot\cdot\mathrm{w(}C_{1})$ and we set
$d(T)=d($\textrm{w(}$T)).$

\begin{example}
$T=%
\begin{tabular}
[c]{|l|lll}\hline
$\mathtt{1}$ & $\mathtt{2}$ & \multicolumn{1}{|l}{$\mathtt{3}$} &
\multicolumn{1}{|l|}{$\mathtt{\bar{1}}$}\\\hline
$\mathtt{4}$ & $\mathtt{4}$ & \multicolumn{1}{|l}{$\mathtt{\bar{3}}$} &
\multicolumn{1}{|l}{}\\\cline{1-3}%
$\mathtt{\bar{4}}$ & $\mathtt{\bar{2}}$ & \multicolumn{1}{|l}{$\mathtt{\bar
{1}}$} & \multicolumn{1}{|l}{}\\\cline{1-3}%
$\mathtt{\bar{3}}$ &  &  & \\\cline{1-1}%
\end{tabular}
$ is a symplectic tableau and \textrm{w}$(T)=\overline{1}3\overline
{3}\,\overline{1}24\overline{2}14\overline{4}\,\overline{3}$
\end{example}

Kashiwara and Nakashima have proved in \cite{KN} the following crucial
theorem, which provides a concrete model for the crystal of an irreducible
$U_{q}(sp_{2n})$-module$.$

\begin{theorem}
\label{thKN}Let $\lambda=\overset{n}{\underset{i=1}{\sum}}\lambda_{i}%
\Lambda_{i}$ be an integral dominant weight. Then the readings of the
symplectic tableaux of shape $\lambda$ are the vertices of a connected
component of the crystal graph $G_{n,\left|  \lambda\right|  }$ isomorphic to
$B(\lambda)$. The highest weight vertex of this graph is the reading of the
tableau of shape $\lambda$ having only $i$'s on its $i$-th row for $i=1,2,...,n$
\end{theorem}

Using Lemma \ref{lem_coplactic}, we deduce immediately the following corollary
which will be used repeatedly.

\begin{corollary}
\label{cortabinword}Let $w,w^{\prime}\in\mathcal{C}_{n}^{\ast}$ such that
$w\longleftrightarrow w^{\prime}$. Write $w=uw_{T}v$ and $w^{\prime}%
=u^{\prime}w_{T}^{\prime}v^{\prime}$ such that $u,u^{\prime},v,v^{\prime
},w_{T},w_{T}^{\prime}\in\mathcal{C}_{n}^{\ast}$, $\frak{L}(u^{\prime
})=\frak{L}(u)$ and $\frak{L}(v^{\prime})=\frak{L}(v)$. Then $w_{T\text{ }}$
is the reading of a symplectic tableau if and only if $w_{T}^{\prime}$ is the
reading of a symplectic tableau.
\end{corollary}

\section{A plactic monoid for $U_{q}(sp_{2n})\label{sec_mono}$}

\subsection{Symplectic tableau associated to a word}

\begin{definition}
Let $w_{1}$ and $w_{2}$ be two words.\ We write $w_{1}\sim w_{2}$ when these
two words occur at the same place in the two isomorphic connected components
$B(w_{1})$ and $B(w_{2})$ of the crystal graph $G_{n}$, that is if there exist
$i_{1},...,i_{r}$ such that $w_{1}=\widetilde{f}_{i_{i}}\cdot\cdot
\cdot\widetilde{f}_{i_{r}}(w_{1}^{0})$ and $w_{2}=\widetilde{f}_{i_{i}}%
\cdot\cdot\cdot\widetilde{f}_{i_{r}}(w_{2}^{0})$, where $w_{1}^{0}$ and
$w_{2}^{0}$ are the vertices of highest weight of $B(w_{1})$ and $B(w_{2})$.
\end{definition}

Then by Theorem \ref{thKN} we have

\begin{proposition}
\label{propunideTab}For any word $w$ of $\mathcal{C}_{n}^{\ast}$ there exists
a unique symplectic tableau $T$ such that $w\sim\mathrm{w(}T).$ We shall
denote it by $T=P(w)$.
\end{proposition}

\subsection{The monoid $Pl(C_{n})$}

\begin{definition}
\label{defmonoidpl}The monoid $Pl(C_{n})$ is the quotient of the free monoid
$\mathcal{C}_{n}^{\ast}$ by the relations:

$R_{1}:$ $yzx\equiv yxz$ for $x\leq y<z$ with $z\neq\overline{x},$ and
$xzy\equiv zxy$ for $x<y\leq z$ with $z\neq\overline{x};$

$R_{2}:y(\overline{x-1})(x-1)\equiv yx\overline{x}$ and $x\overline{x}%
y\equiv(\overline{x-1})(x-1)y$ for $1<x\leq n$ and $x\leq y\leq\overline{x};$

$R_{3}:$ let $w$ be a non admissible column word such that each strict factor
of $w$ is an admissible column word. Write $z$ for the lowest unbarred letter
such that the pair $(z,\overline{z})$ occurs in $w$ and $N(z)=z+1$ (see Remark
\ref{remNz}). Then $w\equiv\widetilde{w}$ where $\widetilde{w}$ is the column
word obtained by erasing the pair $(z,\overline{z})$ in $w.$
\end{definition}

It is clear that $w_{1}\equiv w_{2}$ implies $d(w_{1})=d(w_{2})$, that is,
$\equiv$ is compatible with the grading given by $d$. Notice that $R_{1}$
contains the Knuth relations for type $A$. In $R_{3}$ the condition $N(z)>z$
is equivalent to $N(z)=z+1.$ Indeed it is obvious if $z=1$ and, if $z\neq1,$
$(z,\overline{z})$ is contained in a strict subword of $w$ of length
$\frak{L}(w)-1$ which is admissible.

Denote by $\xi$ the crystal isomorphism $B(121)\overset{\cong}{\rightarrow
}B(112)$. The words occurring in the right hand side of $R_{1}$ and $R_{2}$
are the vertices of $B(112)$. Indeed $112$ is the reading of a symplectic
tableau of shape $\Lambda_{2}+\Lambda_{1}$.\ So by Theorem \ref{thKN},
$B(112)$ contains all the readings of the symplectic tableaux of this shape.
By Theorem \ref{thKN} the connected component of $G_{n}$ containing $21$ is
$B(11)$ which contains the words $bc$ with $b\geq c$. Moreover $B(12)$
contains the words $ab\neq1\overline{1}$ with $a<b$. Hence $B(121)$ contains
words $abc$ with $a<b$, $b\geq c$, and $ab\neq1\overline{1}.$ But $B(121)\cong
B(112)$ so these connected components have the same number of vertices. So
$B(121)$ contains all the words of the previous type, that is, exactly the
words occurring in the left hand side of relations $R_{1}$ and $R_{2}$.

\begin{lemma}
\label{Lem_B(112)_iso_B(121)}Let $w$ be a word occurring in the left hand side
of a relation $R_{1}$ or $R_{2}$. Then $\xi(w)$ is the word occurring in the
right hand side of this relation.
\end{lemma}

\begin{proof}
The lemma is true for $w=121$. Consider $w\in B(112)$ verifying the lemma and
such that $\widetilde{f}_{i}(w)\neq0$.\ By induction, it suffices to prove
that the words $v=\widetilde{f}_{i}(w)$ and $\xi(v)=\widetilde{f}_{i}(\xi(w))$
are respectively the left and right hand sides of the same relation $R_{1}$ or
$R_{2}$. In the sequel we restrict ourselves to the case $w=yzx$ with $x\leq
y<z$.\ The case $w=xzy$ with $x<y\leq z$ may be treated similarly. Since $w$
verifies the lemma, one of the two letters $y$ or $x$ is not modified when we
compute $\xi(w)$ from $w$.\ More precisely, when $w=yzx$ is not of the form
$x\overline{x}x,$ we can write
\[
\xi(w)=yx^{\prime}z^{\prime}\text{ with }\left\{  \text{%
\begin{tabular}
[c]{l}%
$x^{\prime}=x$ and $z^{\prime}=z$ if $z\neq\overline{x}$,\\
$x^{\prime}=(x+1)$ and $z^{\prime}=(\overline{x+1})$ otherwise.
\end{tabular}
}\right.
\]
When $w=x\overline{x}x$, we have $\xi(w)=(\overline{x-1})(x-1)x$ with $1<x\leq
n.$ Since $\widetilde{f}_{i}(w)\neq0,$ $i=x$.\ For $i<n,$ we obtain
$v=x\overline{x}(x+1)$ and $\xi(v)=\widetilde{f}_{i}\xi(w)=(\overline
{x-1})(x-1)(x+1)$.\ For $i=n,$ we have $v=n\overline{n}\,\overline{n}$ and
$\xi(v)=\widetilde{f}_{n}\xi(w)=(\overline{n-1})(n-1)\overline{n}$. In both
cases the lemma is verified.

Now suppose that $w=xyz$ is not of the form $x\overline{x}x.$ Then
$zx\neq\overline{n}n$ (otherwise $w=n\overline{n}n$).\ When we apply
$\widetilde{f}_{i}$ to $w$ one of the three letters $t\in\{y,x,z\}$ is changed
into $\widetilde{f}_{i}(t)$.\ If $v=\widetilde{f}_{i}(y)zx$, then we obtain by
(\ref{tens}) $\varepsilon_{i}(z)=0$ and $\xi(v)=\widetilde{f}_{i}%
\xi(w)=\widetilde{f}_{i}(y)x^{\prime}z^{\prime}$.\ Indeed we have
$\varepsilon_{i}(z^{\prime}x^{\prime})=\varepsilon_{i}(zx)$ and $\varphi
_{i}(z^{\prime}x^{\prime})=\varphi_{i}(zx).$ So $x<\widetilde{f}_{i}(y)<z$ and
$v$ verifies the lemma. If $v=y\widetilde{f}_{i}(zx)$, we have similarly
$\xi(v)=y\widetilde{f}_{i}(x^{\prime}z^{\prime})$ and we consider the two
following cases:

\noindent$\mathrm{(i)}:v=y\widetilde{f}_{i}(z)x$.\ Then $\varepsilon_{i}(x)=0$
hence $z\neq\overline{x}$ and $\xi(w)=yxz$.

\noindent When $\varphi_{i}(x)=0,$ we have $\varepsilon_{i}(\overline{x})=0$
hence $\widetilde{f}_{i}(z)\neq\overline{x}$.\ We obtain $\xi(v)=\widetilde
{f}_{i}\xi(w)=yx\widetilde{f}_{i}(z)$ with $x\leq y<\widetilde{f}_{i}(z)$ and
$\widetilde{f}_{i}(z)\neq\overline{x}$. So $v$ verifies the lemma

\noindent When $\varphi_{i}(x)=1,$ we have $i<n$ because $z>x$ and
$\varphi_{i}(z)=1$. Therefore $x<n$ and $z=\overline{x+1}$.\ So $w=y(\overline
{x+1})x,$ $\xi(w)=y(\overline{x+1})x,$ $v=y\overline{x}x$.\ If $y\neq x$, we
obtain $\xi(v)=\widetilde{f}_{i}\xi(w)=y(x+1)(\overline{x+1)}$.\ If $y=x,$ we
have $\xi(v)=(x+1)(\overline{x+1})x$. In both cases the lemma is verified.

\noindent$\mathrm{(ii)}:v=yz\widetilde{f}_{i}(x)$.\ Then $\varphi_{i}(z)=0$.

\noindent When $z=\overline{x},$ we have $x<n$ (since $zx\neq\overline{n}n$)
and $x<y$ (since $yzx\neq x\overline{x}x$).\ So $i<n$ and $\xi(w)=y(\overline
{x+1})(x+1)$.\ We obtain $v=y\overline{x}(x+1)$ and $\xi(v)=\widetilde{f}%
_{i}\xi(w)=y\overline{x}(x+1)$ with $x+1\leq y<\overline{x}$. So the lemma is verified.

\noindent When $z\neq\overline{x},$ we have $\xi(w)=yxz$. If $x=y,$ we have by
(\ref{tens}) $\varepsilon_{i}(z)=1$.\ So $z=\widetilde{f}_{i}(x)$ (since
$z\neq\overline{x}$).\ We obtain $w=x\widetilde{f}_{i}(x)x$, $\xi
(w)=xx\widetilde{f}_{i}(x),$ $v=x\widetilde{f}_{i}(x)\widetilde{f}(x)$ and
$\xi(v)=\widetilde{f}_{i}(x)x\widetilde{f}_{i}(x)$.\ If $x<y,$ we have
$z\notin\{x+1,\overline{x}\}$ (since $z\geq x+2$) so $\varepsilon_{i}(z)=0.$
We obtain $\xi(v)=y\widetilde{f}_{i}(x)z$ with $\widetilde{f}_{i}(x)\leq y<z$
and $z\neq\overline{\widetilde{f}_{i}(x)}$ (because $\varphi_{i}(z)=0)$.\ In
both cases the lemma is verified.
\end{proof}%

\begin{gather*}%
\begin{tabular}
[c]{lllllllll}%
&  &  &  & $121$ &  &  &  & \\
&  &  & $\overset{\text{1}}{\swarrow}$ &  & $\overset{\text{2}}{\searrow}$ &
&  & \\
&  & $122$ &  &  &  & $1\bar{2}1$ &  & \\
&  & {\tiny 2}$\downarrow$ &  &  &  & {\tiny 1}$\downarrow$ &  & \\
&  & $1\bar{2}2$ &  &  &  & $2\bar{2}1$ &  & \\
& $\overset{\text{2}}{\swarrow}$ & {\tiny 1}$\downarrow$ &  &  &  &
{\tiny 1}$\downarrow$ &  & \\
$1\bar{2}\bar{2}$ &  & $2\bar{2}2$ &  &  &  & $2\bar{1}1$ &  & \\
& $\overset{\text{1}}{\searrow}$ & {\tiny 2}$\downarrow$ &  &  &  &
{\tiny 2}$\downarrow$ & $\overset{\text{1}}{\searrow}$ & \\
&  & $2\bar{2}\bar{2}$ &  &  &  & $\bar{2}\bar{1}1$ &  & $2\bar{1}2$\\
&  & {\tiny 1}$\downarrow$ &  &  &  & {\tiny 1}$\downarrow$ & $\underset
{\text{2}}{\swarrow}$ & \\
&  & $2\bar{1}\bar{2}$ &  &  &  & $\bar{2}\bar{1}2$ &  & \\
&  & {\tiny 1}$\downarrow$ &  &  &  & {\tiny 2}$\downarrow$ &  & \\
&  & $2\bar{1}\bar{1}$ &  &  &  & $\bar{2}\bar{1}\bar{2}$ &  & \\
&  &  & $\overset{\text{2}}{\searrow}$ &  & $\overset{\text{1}}{\swarrow}$ &
&  & \\
&  &  &  & $\bar{2}\bar{1}\bar{1}$ &  &  &  &
\end{tabular}%
\begin{tabular}
[c]{lllllllll}%
&  &  &  & $112$ &  &  &  & \\
&  &  & $\overset{\text{1}}{\swarrow}$ &  & $\overset{\text{2}}{\searrow}$ &
&  & \\
&  & $212$ &  &  &  & $11\bar{2}$ &  & \\
&  & {\tiny 2}$\downarrow$ &  &  &  & {\tiny 1}$\downarrow$ &  & \\
&  & $\bar{2}12$ &  &  &  & $21\bar{2}$ &  & \\
& $\overset{\text{2}}{\swarrow}$ & {\tiny 1}$\downarrow$ &  &  &  &
{\tiny 1}$\downarrow$ &  & \\
$\bar{2}1\bar{2}$ &  & $\bar{1}12$ &  &  &  & $22\bar{2}$ &  & \\
& $\overset{\text{1}}{\searrow}$ & {\tiny 2}$\downarrow$ &  &  &  &
{\tiny 2}$\downarrow$ & $\overset{\text{1}}{\searrow}$ & \\
&  & $\bar{1}1\bar{2}$ &  &  &  & $\bar{2}2\bar{2}$ &  & $22\bar{1}$\\
&  & {\tiny 1}$\downarrow$ &  &  &  & {\tiny 1}$\downarrow$ & $\underset
{\text{2}}{\swarrow}$ & \\
&  & $\bar{1}2\bar{2}$ &  &  &  & $\bar{2}2\bar{1}$ &  & \\
&  & {\tiny 1}$\downarrow$ &  &  &  & {\tiny 2}$\downarrow$ &  & \\
&  & $\bar{1}2\bar{1}$ &  &  &  & $\bar{2}\bar{2}\bar{1}$ &  & \\
&  &  & $\overset{\text{2}}{\searrow}$ &  & $\overset{\text{1}}{\swarrow}$ &
&  & \\
&  &  &  & $\bar{1}\bar{2}\bar{1}$ &  &  &  &
\end{tabular}
\\
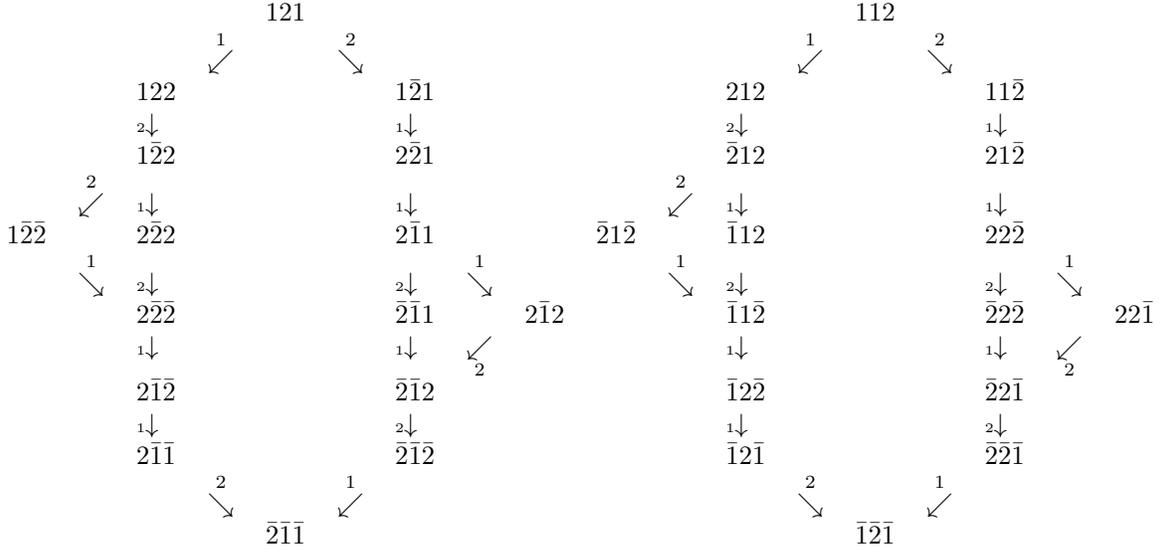
\text{The crystals }B(121)\text{ and }B(112)\text{ for }n=2.
\end{gather*}

If $w=x_{1}\cdot\cdot\cdot x_{p+1}$ is a non admissible column word of length
$p+1$ each strict factor of which is admissible, then $B(w)=B(12\cdot
\cdot\cdot p\overline{p})$.\ Indeed $x_{1}\cdot\cdot\cdot x_{p}$ and
$x_{p}x_{p+1}$ are admissible column words, hence by Lemmas
\ref{lmmaximalweight} and \ref{lem_coplactic} the highest weight vertex of
$B(w)$ is necessarily $12\cdot\cdot\cdot p\overline{p}$.\ Conversely, every
vertex $w\in$ $B(12\cdot\cdot\cdot p\overline{p})$ may be written
$w=z_{1}\cdot\cdot\cdot z_{p}y$ where $z_{1}\cdot\cdot\cdot z_{p}\in
B(1\cdot\cdot\cdot p)$ and $z_{p}y\in B(p\overline{p})$. If $p=1,$
$w=1\overline{1}$ hence we can suppose $p\geq2$.\ Then the highest weight
vertex of $B(p\overline{p})$ is $12$ for $p\overline{p}$ is an admissible
column word of two letters. So we have $z_{1}<\cdot\cdot\cdot<z_{p}<y,$ that
is, $w$ is a column word. The word $w\notin B(1\cdot\cdot\cdot(p+1))$ so it is
not admissible.\ If $v$ is a strict factor of $w$, it is admissible because it
occurs in the same connected component as a strict factor of $12\cdot
\cdot\cdot p\overline{p}$. Finally the words of $B(12\cdot\cdot\cdot
p\overline{p})$ are the non admissible column words of length $p+1$, each
strict factor of which is admissible.

Consider $p>1$. For any word $w\in B(12\cdot\cdot\cdot p\overline{p}),$ we set
$\xi_{p}(w)=\widetilde{w}$ with the notation of Definition \ref{defmonoidpl}.

\begin{lemma}
\label{Lem_iso_B(11bar)_(B0)}The map $\xi_{p}$ is a crystal isomorphism from
$B(12\cdot\cdot\cdot p\overline{p})$ to $B(12\cdot\cdot\cdot(p-1))$.
\end{lemma}

\begin{proof}
We know that $\xi_{p}(12\cdot\cdot\cdot p\overline{p})=12\cdot\cdot\cdot
(p-1)$.\ Hence it suffices to prove that for any $i=1,...,n$ and any $w\in
B(12\cdot\cdot\cdot p\overline{p}),$ such that $\widetilde{f}_{i}(w)\neq0$,
$\xi_{p}\widetilde{f}_{i}(w)=\widetilde{f}_{i}\xi_{p}(w)$. Denote by
$(z,\overline{z})$ the pair of letters erased in $w$ when we compute $\xi
_{p}(w)$.

Suppose $i\notin\{z-1,z\}.$ The words $w$ and $\xi_{p}(w)$ differ only by the
letters $z$ and $\overline{z}$ which do not interfere in the computation of
$\widetilde{f}_{i}$.\ So it suffices to show that the pair $(z,\overline{z})$
disappears in $\widetilde{f}_{i}(w)$ when the relation $R_{3}$ is applied.\ We
have again $N(z)=z+1$ in $\widetilde{f}_{i}(w)$ because $i\notin\{z-1,z\}$.
Suppose that $\widetilde{f}_{i}(w)$ contains a pair $(t,\overline{t})$
satisfying $N(t)=t+1$ and $t<z$. The word $w$ can not contain the pair
$(t,\overline{t})$ otherwise $(z,\overline{z})$ is not the pair of letters
which disappears when we apply $R_{3}$. Hence a letter $t-1\in w$ is replaced
by $t$ when we compute $\widetilde{f}_{i}(w)$. Indeed, by (\ref{regle+-}) the
case $t\in w$ and a letter $\overline{t+1}$ is replaced by $\overline{t}$ is
impossible. Therefore we have $N(t-1)=t+1$ in $w$.$\;$So $w$ contains at least
a strict factor $v$ such that $N(t-1)=t$ in $v$.$\;$Such a factor is not an
admissible column word. So we obtain a contradiction. Hence the pair
$(z,\overline{z})$ disappears in $\widetilde{f}_{i}(w)$ when we compute
$\xi_{p}(w)$.

Suppose $i=z-1$. Then $z\neq1$ and $w$ must contain a letter $z-1$ or
$\overline{z-1}$. Otherwise, if $x$ is the greatest unbarred letter $<z$ such
that $x\in w$ or $\overline{x}\in w$, $N(x)=z-1>x+1$ and the pair
$(z,\overline{z})$ does not disappear by applying $R_{3}$. There are three
cases to consider: $(\mathrm{i})$ $\overline{z-1}$ occurs alone in $w,$ or
$(\mathrm{ii})$ $z-1$ occurs alone in $w,$ or $(\mathrm{iii})$ the pair
$(z-1,\overline{z-1})$ occurs in $w.$ In each case the result is verified by a
direct computation. For example, in case $(\mathrm{ii}),$ $w=\cdot\cdot
\cdot(z-1)z\cdot\cdot\cdot\overline{z}\cdot\cdot\cdot$ we have $\xi
_{p}(w)=\cdot\cdot\cdot(z-1)\cdot\cdot\cdot\cdot\cdot\cdot$ (we have only
written the letters $\overline{z-1},\overline{z},z,z-1$ of $w$ and $\xi
_{p}(w)$). So we obtain $\widetilde{f}_{i}(w)=\cdot\cdot\cdot(z-1)z\cdot
\cdot\cdot\overline{z-1}\cdot\cdot\cdot$. Hence $\xi_{p}\widetilde{f}%
_{i}(w)=\cdot\cdot\cdot z\cdot\cdot\cdot\cdot\cdot\cdot=\widetilde{f}_{i}%
\xi_{p}(w).$

Suppose $i=z$. The letters $z+1$ and $\overline{z+1}$ do not occur
simultaneously in $w$ (otherwise the column word obtained by erasing the last
letter of $w$ is not admissible because $N(z+1)>z+1$ in this word).
Considering the cases: $(\mathrm{i})$ $z+1$ occurs in $w$, $(\mathrm{ii})$
$\overline{z+1}$ occurs in $w$, $(\mathrm{iii})$ neither $z+1$ or
$\overline{z+1}$ occurs in $w,$ we obtain the equality $\xi_{p}\widetilde
{f}_{i}(w)=\widetilde{f}_{i}\xi_{p}(w)$ by a direct computation.
\end{proof}

The following proposition shows the compatibility of the relations $R_{i}$
above with crystal graphs operators:

\begin{proposition}
\label{prop_cop_fi_cong}Let $w_{1}$ and $w_{2}$ be words of $\mathcal{C}%
_{n}^{\ast}$ such that $w_{1}\equiv w_{2}$. Then for $i=1,...,n:$%
\[
\widetilde{e}_{i}(w_{1})\equiv\widetilde{e}_{i}(w_{2})\ ,\varepsilon_{i}%
(w_{1})=\varepsilon_{i}(w_{2})\ ,\widetilde{f}_{i}(w_{1})\equiv\widetilde
{f}_{i}(w_{2})\ ,\varphi_{i}(w_{1})=\varphi_{i}(w_{2}).
\]
\end{proposition}

\begin{proof}
By induction we can suppose that $w_{2}$ is obtained from $w_{1}$ by applying
only one plactic relation. In this case we write $w_{1}=u\widehat{w}_{1}v$ and
$w_{2}=u\widehat{w}_{2}v$ where $\widehat{w}_{1},\widehat{w}_{2}$ are minimal
factors of $w_{1}\ $and $w_{2}$ such that $\widehat{w}_{1}\equiv\widehat
{w}_{2}$. Formula (\ref{tens}) implies that it is enough to prove the
proposition for $\widehat{w}_{1}$ and $\widehat{w}_{2}$. When $\widehat{w}%
_{1}$ and $\widehat{w}_{2}$ differ by applying a relation $R_{1}$ or $R_{2}$
the proposition follows from Lemma \ref{Lem_B(112)_iso_B(121)}.\ When
$\widehat{w}_{1}\in B(1\cdot\cdot\cdot p\overline{p})$ it is a consequence of
Lemma \ref{Lem_iso_B(11bar)_(B0)}.
\end{proof}

\begin{corollary}
\label{cor1}If $w_{1}\equiv w_{2}$ then $w_{1}\sim w_{2}$ (two congruent words
occur at the same place in two isomorphic connected components of $G_{n}$).
\end{corollary}

\begin{proof}
Let $w_{1}^{0}$ and $w_{2}^{0}$ be the highest weight vertices of $B(w_{1})$
and $B(w_{2})$. Then the above proposition and an induction prove that
$w_{1}^{0}\equiv w_{2}^{0}$. So $d(w_{1}^{0})=d(w_{2}^{0})$ hence $B(w_{1})$
and $B(w_{2})$ are isomorphic. Moreover, by the previous proposition there
exists $i_{1},...,i_{r}\in\{1,...,n\}$ such that $w_{1}^{0}=\widetilde
{e}_{i_{1}}\widetilde{e}_{i_{2}}\cdot\cdot\cdot\widetilde{e}_{i_{r}}(w_{1})$
and $w_{2}^{0}=\widetilde{e}_{i_{1}}\widetilde{e}_{i_{2}}\cdot\cdot
\cdot\widetilde{e}_{i_{r}}(w_{2})$.\ So $w_{1}$ and $w_{2}$ occur at the same
place in two isomorphic connected components.
\end{proof}

We shall now prove the converse%
\begin{equation}
w_{1}\sim w_{2}\Longrightarrow w_{1}\equiv w_{2}. \label{(E)}%
\end{equation}

\begin{proposition}
\label{prop_HP_congru}Consider $w$ a highest weight vertex of $G_{n}$. Then
$w\equiv\mathrm{w(}P(w))$ (every highest weight vertex is congruent to a
highest weight tableau word).
\end{proposition}

\begin{proof}
We proceed by induction on $\frak{L}(w)$. If $\frak{L}(w)=1$ then $w=1$ and
$P(w)=%
\begin{tabular}
[c]{|l|}\hline
$1$\\\hline
\end{tabular}
$. Suppose the proposition true for the highest weight vertices of length
$\leq l$ and consider a highest weight vertex $w$ of length $l+1$. Write
$w=vx$ where $x\in\mathcal{C}_{n}$ and $v$ is a word of length $l$. Then by
Lemma \ref{lmmaximalweight} $v$ is of highest weight so by induction
$v\equiv\mathrm{w(}P(v))$ where $P(v)$ contains only unbarred letters. Hence
$w\equiv\mathrm{w(}P(v))x$. We deduce from the condition $\varepsilon
_{i}(x)\leq\varphi_{i}(v)$ of Lemma \ref{lmmaximalweight} that only the two
following situations can occur:

\begin{enumerate}
\item $x$ is an unbarred letter and $P(v)$ contains a column of reading
$1\cdot\cdot\cdot(x-1)$ if $x>1,$

\item $x=\overline{p}$ and $P(v)=C_{1}\cdot\cdot\cdot C_{r}$ contains a column
$C_{i}$ with $\mathrm{w(}C_{i})=1\cdot\cdot\cdot p$ (we suppose $i$ maximal in
the sequel).
\end{enumerate}

$1:$ If $x$ is unbarred $P(w)$ is obtained from $P(v)$ by adding a box
containing $x$ at the bottom of the leftmost column of reading $1\cdot
\cdot\cdot(x-1)$ if $x>1,$ by adding a box containing $1$ to the right of
$P(v)$ otherwise. We have $w\equiv vx\equiv\mathrm{w(}P(v))x\equiv
\mathrm{w(}P(w))$ because in $Pl(C_{n})$, $x$ commute with all the column
words containing $x$ and only unbarred letters.

$2:$ If $x=\overline{p}$, $P(w)$ is the reading of the tableau obtained by
erasing the letter $p$ on the top of $C_{i}$.\ The word \textrm{w}$(C_{1})$ is
of highest weight, so we can write \textrm{w}$(C_{1})=1\cdot\cdot\cdot q$ with
$q\geq p$ and $vx=$\textrm{w}$(C_{r})\cdot\cdot\cdot$\textrm{w}$(C_{2}%
)(1\cdot\cdot\cdot q\overline{p})$. Then by using the contraction relation
$1\cdot\cdot\cdot q\overline{p}\equiv1\cdot\cdot\cdot\widehat{p}\cdot
\cdot\cdot q$ (where the hat means removal the letter $p$) we obtain $w\equiv
$\textrm{w}$(C_{r})\cdot\cdot\cdot$\textrm{w}$(C_{2})1\cdot\cdot\cdot
\widehat{p}\cdot\cdot\cdot q$. We have $P(w)=P($\textrm{w}$(C_{r})\cdot
\cdot\cdot$\textrm{w}$(C_{2})1\cdot\cdot\cdot\widehat{p}\cdot\cdot\cdot
q)$.\ Hence by case $1$ and an easy induction $w\equiv\mathrm{w(}P(w))$.
\end{proof}

\begin{corollary}
If $w_{1}^{0}$ and $w_{2}^{0}$ are two highest weight vertices with the same
weight, then $w_{1}^{0}\equiv w_{2}^{0}$
\end{corollary}

\begin{proof}
$P(w_{1}^{0})=P(w_{2}^{0})=T$ because $w_{1}^{0}$ and $w_{2}^{0}$ have the
same weight. Then $w_{1}^{0}\equiv T\equiv w_{2}^{0}$.
\end{proof}

\begin{theorem}
\label{Th_good_rela}For any vertices $w_{1}$ and $w_{2}$ of $G_{n}$%
\[
w_{1}\sim w_{2}\Longleftrightarrow w_{1}\equiv w_{2}.
\]
\end{theorem}

\begin{proof}
Suppose that $w_{1}\sim w_{2}$ and denote by $w_{1}^{0}$ and $w_{2}^{0}$ the
highest weight vertices of $B(w_{1})$ and $B(w_{2}).$ The corollary above
shows that $w_{1}^{0}\equiv w_{2}^{0}$. Write $w_{1}=\widetilde{f}_{i_{1}%
}\cdot\cdot\cdot\widetilde{f}_{i_{r}}w_{1}^{0}$, then we have $w_{2}%
=\widetilde{f}_{i_{1}}\cdot\cdot\cdot\widetilde{f}_{i_{r}}w_{2}^{0}$. Hence by
Proposition \ref{prop_cop_fi_cong}, $w_{1}\equiv w_{2}$. The converse was
proved in Corollary \ref{cor1}.
\end{proof}

\noindent\textbf{Remark }The sufficiency of the relations obtained from the
crystal isomorphisms $B(121)\cong B(112)$ and $B(12\cdot\cdot\cdot
p\overline{p})\overset{\cong}{\rightarrow}B(12\cdot\cdot\cdot p)$ to generate
$\mathcal{C}_{n}^{\ast}/\sim$ was proved in \cite{Lit}.

\section{A bumping algorithm for type $C$}

Now we are going to see how the symplectic tableau $P(w)$ may be computed for
each vertex $w$ by using an insertion scheme analogous to the bumping
algorithm for type $A$. As a first step, we compute $P(w)$ when $w=$%
\textrm{w(}$C)x$, where $x$ and $C$ are respectively a letter and an
admissible column. This will be called ``the insertion of the letter $x$ in
the admissible column $C$'' and denoted by $x\rightarrow C=P(w)$.

Then we will be able to obtain $P(w)$ when $w=$\textrm{w(}$T)x$ with $x$ a
letter and $T$ a symplectic tableau. This will be called ``the insertion of
the letter $x$ in the symplectic tableau $T$'' and denoted by $x\rightarrow
T$. Our construction of $P$ will be recursive, in the sense that if $P(u)=T$
and $x\in\mathcal{C}_{n}$, then $P(ux)=x\rightarrow T$.

\subsection{Insertion of a letter in an admissible column\label{xC}}

Consider a word $w=$\textrm{w(}$C)x$, where $x$ and $C$ are respectively a
letter and an admissible column of height $p$. Then by Lemma
\ref{lmmaximalweight} the highest weight vertex $w^{0}$ of $B(w)\,\ $verifies:%
\begin{align*}
(\text{\textrm{i)}}  &  :\text{ }w^{0}=1\cdot\cdot\cdot p\,(p+1)\text{ or}\\
(\mathrm{ii)}  &  :\text{ }w^{0}=w^{0}=1\cdot\cdot\cdot p\,\overline{p}\text{
or }\\
(\mathrm{iii)}  &  :\text{ }w^{0}=1\cdot\cdot\cdot p\,1
\end{align*}
In case $(\mathrm{i}$\textrm{) }$w^{0}$ is an admissible column word so $w$ is
the reading of the admissible column obtained by adding a box filled by $x$ at
the bottom of $C$ that is $x\rightarrow C=%
\begin{tabular}
[c]{|l|}\hline
$C$\\\hline
$x$\\\hline
\end{tabular}
$.

\noindent In case $(\mathrm{ii})$\textrm{ }we have seen in the previous
section that the vertex $w$ is a non admissible column word all of whose
proper factors are admissible. Hence $x\rightarrow C$ is the column of reading
$\widetilde{w}$ obtained from $w$ by a congruence $R_{3}$ as described in
Definition \ref{defmonoidpl}.\ 

\noindent In case $(\mathrm{iii),}$ the congruence class of $w^{0}$ consists
of the words $w_{i}^{0}=12\cdot\cdot\cdot i1(i+1)(i+2)\cdot\cdot\cdot p$ for
$p\geq i>0,$ and the congruence relations on this class define a linear graph
with end points $w^{0}=w_{p}^{0}$ and $w_{1}^{0}=112\cdot\cdot\cdot p$ (it
describes the commutation of the letter $1$ with the column word $12\cdot
\cdot\cdot p;$ with all congruences of type $R_{1}$ or $R_{2}$).\ By repeated
applications of operators $\widetilde{f}_{i},$ this class can be transformed
into that of $w,$ which will also have a linear graph with all congruences of
type $R_{1}$ or $R_{2}$, and with a word $w^{\prime}\in B(w_{1}^{0})$ at the
end opposite to $w$.\ Then $w^{\prime}$ is the reading of a symplectic tableau
consisting of a column $C^{\prime}$ of height $p$ and a column
\begin{tabular}
[c]{|l|}\hline
$x^{\prime}$\\\hline
\end{tabular}
(with $x^{\prime}\in\mathcal{C}_{n}$). We can write $x\rightarrow C=C^{\prime
}$%
\begin{tabular}
[c]{|l|}\hline
$x^{\prime}$\\\hline
\end{tabular}
$=P(w)$.\ To obtain the congruence class of $w,$ one applies a congruence of
type $R_{1}$ or $R_{2}$ to the final sub-word of length $3$ of $w$ (there is
only one possibility and it applies from left to right). On the resulting
word, one continues with the overlapping sub-word of length $3$ one place to
the left and so forth until the left subword of length $3$ has been operate
upon (see Example \ref{ex_x_inC} below). This insertion can be regarded as the
combinatorial description of the crystal isomorphism:%
\begin{equation}
B(1\cdot\cdot\cdot p1)\rightarrow B(11\cdot\cdot\cdot p)
\label{iso_letter_col}%
\end{equation}

\begin{example}
\label{ex_x_inC}Suppose $w=35\bar{5}\bar{4}\bar{3}.$

If $x=\bar{2}$ the word $35\bar{5}\bar{4}\bar{3}\bar{2}$ is a non admissible
column word each strict subword of which is an admissible column word. Then we
obtain by applying $R_{3}$, $\bar{2}\rightarrow%
\begin{tabular}
[c]{|l|}\hline
$\mathtt{3}$\\\hline
$\mathtt{5}$\\\hline
$\mathtt{\bar{5}}$\\\hline
$\mathtt{\bar{4}}$\\\hline
$\mathtt{\bar{3}}$\\\hline
\end{tabular}
=%
\begin{tabular}
[c]{|l|}\hline
$\mathtt{3}$\\\hline
$\mathtt{\bar{4}}$\\\hline
$\mathtt{\bar{3}}$\\\hline
$\mathtt{\bar{2}}$\\\hline
\end{tabular}
$.

If $x=3$ we obtain by applying at each step $R_{1}$ or $R_{2}$ the new word:%
\[
35\bar{5}\mathbf{\bar{4}\bar{3}3}\equiv35\mathbf{\bar{5}\bar{4}4}\bar{4}%
\equiv3\mathbf{5\bar{5}5}\bar{5}\bar{4}\equiv\mathbf{3\bar{4}4}5\bar{5}\bar
{4}\equiv\bar{4}\,345\bar{5}\bar{4}%
\]
This last word is the reading of a symplectic tableau. \label{exinsert-+}Hence
$3\rightarrow%
\begin{tabular}
[c]{|l|}\hline
$\mathtt{3}$\\\hline
$\mathtt{5}$\\\hline
$\mathtt{\bar{5}}$\\\hline
$\mathtt{\bar{4}}$\\\hline
$\mathtt{\bar{3}}$\\\hline
\end{tabular}
=%
\begin{tabular}
[c]{|l|l}\hline
$\mathtt{3}$ & \multicolumn{1}{|l|}{$\mathtt{\bar{4}}$}\\\hline
$\mathtt{4}$ & \\\cline{1-1}%
$\mathtt{5}$ & \\\cline{1-1}%
$\mathtt{\bar{5}}$ & \\\cline{1-1}%
$\mathtt{\bar{4}}$ & \\\cline{1-1}%
\end{tabular}
$.
\end{example}

\subsection{Insertion of a letter in a symplectic tableau\label{xT}}

Let $T=C_{1}\cdot\cdot\cdot C_{r}$ be a symplectic tableau with admissible
columns $C_{i}$ and $x$ a letter. We have seen in the proof of Proposition
\ref{prop_HP_congru} that the vertex of highest weight of the connected
component containing \textrm{w(}$T)x$ may be written $w^{0}x^{0}$ where
$w^{0}$ is the reading of a highest weight tableau $T^{0}$, and that one of
the following conditions is satisfied:%
\begin{align*}
(\mathrm{i)}  &  :x^{0}=p\text{ and }1\cdot\cdot\cdot(p-1)\text{ is the
reading of a column of }T^{0}\\
(\mathrm{ii)}  &  :x^{0}=\overline{p}\text{ and }1\cdot\cdot\cdot p\text{ is
the reading of a column of }T^{0}\text{.}%
\end{align*}
Write $T^{0}=C_{1}^{0}\cdot\cdot\cdot C_{r}^{0}$. In case ($\mathrm{i)}$ let
$k$ be minimal such that \textrm{w}($C_{k}^{0})=1\cdot\cdot\cdot(p-1)$. Denote
by $T^{0^{\prime}}$ the tableau obtained from $T^{0}$ by adding a box
containing the letter $p$ on the top of $C_{k}^{0}$.\ There exists a unique
sequence $w_{0},...,w_{k-1}$ such that $\mathrm{w}(T^{0^{\prime}})$ is the
highest weight vertex of $B(w_{k-1})$ and such that $w_{j-1}=\mathrm{w}%
(C_{r})\cdot\cdot\cdot\mathrm{w}(C_{j})x_{j-1}\mathrm{w}(C_{j-1}^{\prime
})\cdot\cdot\cdot\mathrm{w}(C_{1}^{\prime})$ is transformed into the congruent
word $w_{j}=\mathrm{w}(C_{r})\cdot\cdot\cdot\mathrm{w}(C_{j+1})x_{j}%
\mathrm{w}(C_{j}^{\prime})\cdot\cdot\cdot\mathrm{w}(C_{1}^{\prime})$ where
$x_{j}$ and $C_{j}^{\prime}$ are determined by $x_{j-1}\rightarrow C_{j}%
=C_{j}^{\prime}%
\begin{tabular}
[c]{|l|}\hline
$x_{j}$\\\hline
\end{tabular}
$.\ The word $w_{k-1}$ is the reading of a symplectic tableau $T^{\prime}$ and
is obtained by computing $k-1$ insertions of type ``insertion of a letter in
an admissible column''. So we have $x\rightarrow T=T^{\prime}$.\ Notice that
there is no contraction in this case because no relation of type $R_{3}$ is applied.

In case ($\mathrm{ii)}$\textrm{, }$\mathrm{w}(C_{1})x$ is a non admissible
column word whose strict subwords are admissible column words. Suppose
$\widetilde{\mathrm{w}(C_{1})x}=y_{1}\cdot\cdot\cdot y_{s}$ and write
$\widehat{T}$\ for the tableau whose columns are $C_{2},...,C_{r}$. Then we
have%
\[
x\rightarrow T=y_{s}\rightarrow(y_{s-1}\rightarrow(\cdot\cdot\cdot
y_{1}\rightarrow\widehat{T}))
\]
that is $x\rightarrow T$ is obtained by inserting successively the letters of
$\widetilde{\mathrm{w}(C_{1})x}$ into the tableau $\widehat{T}$. Indeed there
is no contraction during the insertion of the letters $y_{i}$ because the
highest weight tableau associated to $x\rightarrow T$ is obtained by inserting
the letters $1,...,p-1,p+1,...,s+1$ in the highest weight tableau $C_{2}%
^{0}\cdot\cdot\cdot C_{r}^{0}$ and no relation $R_{3}$ is required to realize
it. So $y_{s}\rightarrow(y_{s-1}\rightarrow(\cdot\cdot\cdot y_{1}%
\rightarrow\widehat{T}))$ is well defined and is a symplectic tableau whose
reading is congruent to \textrm{w(}$T)x$. It must be equal to $x\rightarrow T$.

\begin{example}
let $T=%
\begin{tabular}
[c]{|l|l|ll}\hline
$\mathtt{1}$ & $\mathtt{3}$ & $\mathtt{\bar{3}}$ &
\multicolumn{1}{|l|}{$\mathtt{\bar{2}}$}\\\hline
$\mathtt{\bar{3}}$ & $\mathtt{\bar{3}}$ & $\mathtt{\bar{2}}$ &
\multicolumn{1}{|l}{}\\\cline{1-3}\cline{2-3}%
$\mathtt{\bar{2}}$ & $\mathtt{\bar{1}}$ &  & \\\cline{1-2}\cline{2-2}%
\end{tabular}
$and $x=\overline{1}$. The column $C_{1}=%
\begin{tabular}
[c]{|l|}\hline
$\mathtt{1}$\\\hline
$\mathtt{\bar{3}}$\\\hline
$\mathtt{\bar{2}}$\\\hline
$\mathtt{\bar{1}}$\\\hline
\end{tabular}
$ is not admissible and $\widetilde{C}_{1}=%
\begin{tabular}
[c]{|l|}\hline
$\mathtt{\bar{3}}$\\\hline
$\mathtt{\bar{2}}$\\\hline
\end{tabular}
$.\ So we have to insert $\overline{3}$ next $\overline{2}$ in the tableau
$\widehat{T}=%
\begin{tabular}
[c]{|l|ll}\hline
$\mathtt{3}$ & $\mathtt{\bar{3}}$ & \multicolumn{1}{|l|}{$\mathtt{\bar{2}}$%
}\\\hline
$\mathtt{\bar{3}}$ & $\mathtt{\bar{2}}$ & \multicolumn{1}{|l}{}\\\cline{1-2}%
$\mathtt{\bar{1}}$ &  & \\\cline{1-1}%
\end{tabular}
$. An easy computation shows that the insertion of $\overline{3}$ gives the
tableau
\begin{tabular}
[c]{|l|lll}\hline
$\mathtt{2}$ & $\mathtt{\bar{3}}$ & \multicolumn{1}{|l}{$\mathtt{\bar{2}}$} &
\multicolumn{1}{|l|}{$\mathtt{\bar{2}}$}\\\hline
$\mathtt{\bar{3}}$ & $\mathtt{\bar{2}}$ & \multicolumn{1}{|l}{} &
\\\cline{1-2}%
$\mathtt{\bar{1}}$ &  &  & \\\cline{1-1}%
\end{tabular}
then by inserting $\overline{2}$ we obtain $\overline{1}\rightarrow T=%
\begin{tabular}
[c]{|l|l|ll}\hline
$\mathtt{2}$ & $\mathtt{\bar{3}}$ & $\mathtt{\bar{2}}$ &
\multicolumn{1}{|l|}{$\mathtt{\bar{2}}$}\\\hline
$\mathtt{\bar{3}}$ & $\mathtt{\bar{2}}$ &  & \\\cline{1-2}\cline{2-2}%
$\mathtt{\bar{2}}$ & $\mathtt{\bar{1}}$ &  & \\\cline{1-2}%
\end{tabular}
$.
\end{example}

Note that the insertion scheme of the letter $x$ in the symplectic tableau
$T=C_{1}\cdot\cdot\cdot C_{r}$ obtained above and that described by Baker in
\cite{ba} coincide if no contraction appears during $x\rightarrow C_{1}%
$.\ When such a contraction appears, Baker computes first a admissible skew
tableau $T^{\prime}$ (see Definition \ref{def_skew_tab}) from $T$ by deleting
the box lying on the top of $C_{1}$ and by filling the column diagram so
obtained with the letters of $\widetilde{\mathrm{w}(C_{1})x}$. For example,
with the tableau $T$ above and $x=\overline{1}$ we obtain
\[
T^{\prime}=%
\begin{tabular}
[c]{l|l|ll}\cline{2-4}%
& $\mathtt{3}$ & $\mathtt{\bar{3}}$ & \multicolumn{1}{|l|}{$\mathtt{\bar{2}}$%
}\\\hline
\multicolumn{1}{|l|}{$\mathtt{\bar{3}}$} & $\mathtt{\bar{3}}$ & $\mathtt{\bar
{2}}$ & \multicolumn{1}{|l}{}\\\cline{1-2}\cline{1-3}%
\multicolumn{1}{|l|}{$\mathtt{\bar{2}}$} & $\mathtt{\bar{1}}$ &  &
\\\cline{1-2}%
\end{tabular}
\]
Then the rest of his algorithm may be regarded as a special case of the
Symplectic Jeu de Taquin described in Section \ref{Sec_JDT}.

\subsection{Computation of $P(w)$ for any word $w$}

Consider a word $w\in\mathcal{C}_{n}^{\ast}$. Then we have%
\begin{align*}
P(w)  &  =%
\begin{tabular}
[c]{|l|}\hline
$w$\\\hline
\end{tabular}
\text{ if }w\text{ is a letter,}\\
P(w)  &  =x\rightarrow P(u)\text{ if }w=ux\text{ with }u\text{ a word and
}x\text{ a letter.}%
\end{align*}
that is $P(w)$ may be computed by induction by using the insertion schemes of
Subsections \ref{xC} and \ref{xT}. Notice that when $w=x_{1}\cdot\cdot\cdot
x_{l}$ is the reading of the symplectic tableau $T$, $P(w)=T$ but the equality%
\[
T=x_{l}\rightarrow(x_{l-1}\rightarrow(\cdot\cdot\cdot\rightarrow%
\begin{tabular}
[c]{|l|}\hline
$x_{1}$\\\hline
\end{tabular}
))
\]
is not combinatorially trivial.

\section{Robinson-Schensted type correspondence for type $C_{n}$}

In this section a bijection is established between words $w$ of length $l$ on
$\mathcal{C}_{n}$ and pairs $(P(w),Q(w))$ where $P(w)$ is the symplectic
tableau defined in section 2 and $Q(w)$ is an oscillating tableau. Such a
one-to-one correspondence has already been obtained by Berele \cite{Be} using
King's definition of symplectic tableaux \cite{King} and an appropriate
insertion algorithm. Unfortunately we do not know if this correspondence is
compatible with a monoid structure. Our bijection based on the previous
insertion algorithm will be different from Berele's one but compatible with
the plactic relations $R_{i}$.

We now recall the definition of an oscillating tableau and we give the
construction of $Q(w).$ Next the Robinson-Schensted type correspondence
theorem will be proved.

\subsection{Oscillating tableau associated to a word $w$}

\begin{definition}
(Berele) An oscillating tableau $Q$ of length $l$ is a sequence of Young
diagrams of shape $(Q_{1},...,Q_{l})$ such that any two consecutive shapes
differ by exactly one box $($i.e. $Q_{k+1}/Q_{k}=%
\begin{tabular}
[c]{|l|}\hline
\\\hline
\end{tabular}
$or $Q_{k}/Q_{k+1}=%
\begin{tabular}
[c]{|l|}\hline
\\\hline
\end{tabular}
)$.
\end{definition}

If, for any $i=1,...,l,$ no column of $Q_{i}$ has height greater than $n$ we
will say that $Q$ is an $n$-oscillating tableau. Let $w=x_{1}\cdot\cdot\cdot
x_{l}\in\mathcal{C}_{n}^{\ast}$. The construction of $P(w)$ involves the
construction of the $l$ symplectic tableaux $(P_{1},..,P_{l})$ defined by
$P_{i}=P(x_{l-i+1}\cdot\cdot\cdot x_{l}).$

\begin{notation}
.We denote by $Q(w)$ the sequence of shapes of the symplectic tableaux
$(P_{1},...,P_{l})$ obtained during the construction of $P(w).$
\end{notation}

\begin{proposition}
For any word $w\in\mathcal{C}_{n}^{\ast},$ $Q(w)$ is an $n$-oscillating tableau.
\end{proposition}

\begin{proof}
Write $w=vx$ where $v\in\mathcal{C}_{n}^{\ast}$ and $x\in\mathcal{C}_{n}$. By
induction it suffices to show that the shapes of $P(w)$ and $P(v)$ differ by
exactly one box.\ The shape of $P(w)$ depends only of the connected component
containing $w$. So we can suppose $w$ of highest weight. Then by Lemma
\ref{lmmaximalweight}, $v$ is of highest weight.

Suppose $x=\overline{i}\geq\overline{n},$ then its weight is equal to
$\Lambda_{i-1}-\Lambda_{i}$ and if we denote by $(\lambda_{1},...,\lambda
_{n})$ the coordinates of the weight of $v$ on the basis of fundamental
weights we have $\lambda_{i}>0$ (because $w$ must be of highest weight). By
Theorem \ref{thKN}, we can deduce that during the insertion $x\rightarrow
P(v)$ a column of height $i$ (corresponding to the weight $\Lambda_{i}$) is
turned into a column of height $i-1$ (corresponding to the weight
$\Lambda_{i-1}$). So the shape of $P(w)$ is obtained by erasing one box from
the shape of $T$.

When $x=i\leq n,$ its weight is equal to $\Lambda_{i}-\Lambda_{i-1}$ and
similar arguments show that the shape of $P(w)$ is obtained by adding one box
to the shape of $P(v)$.
\end{proof}

\subsection{Correspondence theorem}

\begin{proposition}
\label{propmccmQ}Let $w_{1}$ and $w_{2}$ be two words of $\mathcal{C}%
_{n}^{\ast}.$ Then
\[
w_{1}\longleftrightarrow w_{2}\Longleftrightarrow Q(w_{1})=Q(w_{2}).
\]
\end{proposition}

\begin{proof}
We proceed by induction on the length $l$ of the words $w_{1}$ and $w_{2}$. If
$l=1$ the result is immediate. If $w_{1}$ and $w_{2}$ have length $l\geq1,$ we
can write $w_{1}=u_{1}x_{1}$ and $w_{2}=u_{2}x_{2}$ with $x_{1},x_{2}$ letters
and $u_{1},u_{2}$ words of length $l-1$. Let $w_{1}^{0}=u_{1}^{0}x_{1}^{0}$
and $w_{2}^{0}=u_{2}^{0}x_{2}^{0}$ be the highest weight vertices of
$B(w_{1})$ and $B(w_{2})$. Write $Q_{1}$ and $Q_{2}$ for the shapes of
$P(w_{1})$ and $P(w_{2})$ (that is those of $P(w_{1}^{0})$ and $P(w_{2}^{0}%
)$).\ We suppose the Proposition true for the words of length $l-1.$ First we
have:%
\[
w_{1}\longleftrightarrow w_{2}\Longleftrightarrow\left\{
\begin{tabular}
[c]{l}%
$u_{1}\longleftrightarrow u_{2}$\\
$Q_{1}=Q_{2}$%
\end{tabular}
\right.
\]
Indeed if $w_{1}\longleftrightarrow w_{2}$ then $u_{1}\longleftrightarrow
u_{2}$ follows from Lemma \ref{lem_coplactic} and we obtain $Q_{1}=Q_{2}$
because the readings of $P(w_{1})$ and $P(w_{2})$ are in the same connected
component of $G_{n}$. Conversely, $u_{1}\longleftrightarrow u_{2}$ implies
that $u_{1}^{0}=u_{2}^{0}$ and it follows from the equality $Q_{1}=Q_{2}$ that
$\mathrm{wt}(w_{1}^{0})=\mathrm{wt}(w_{2}^{0})$ (because the shape of
$P(w_{i}^{0})$ $i=1,2$ coincides with the weight $\mathrm{wt}(w_{i}^{0})$)$.$
So $x_{1}^{0}=x_{2}^{0}$. This means that $w_{1}^{0}=w_{2}^{0}$ i.e.
$w_{1}\longleftrightarrow w_{2}.$ Finally we obtain by induction:%
\[
w_{1}\longleftrightarrow w_{2}\Longleftrightarrow\left\{
\begin{tabular}
[c]{l}%
$Q(u_{1})=Q(u_{2})$\\
$Q_{1}=Q_{2}$%
\end{tabular}
\right.  \Longleftrightarrow Q(w_{1})=Q(w_{2}).
\]
\end{proof}

The Robinson-Schensted correspondence for type $C_{n}$ follows immediately:

\begin{theorem}
\label{TH_RS}Let $\mathcal{C}_{n,l}^{\ast}$ and $\mathcal{O}_{l}$ be the set
of words of length $l$ on $\mathcal{C}_{n}$ and the set of pairs $(P,Q)$ where
$P$ is a symplectic tableau and $Q$ an $n$-oscillating tableau of length $l$
such that $P$ has shape $Q_{l}$ ($Q_{l}$ is the last shape of $Q$). Then the
map:%
\begin{align*}
\Psi &  :\mathcal{C}_{n,l}^{\ast}\rightarrow\mathcal{O}_{l}\\
w &  \mapsto(P(w),Q(w))
\end{align*}
is a bijection.
\end{theorem}

\begin{proof}
By Theorem \ref{Th_good_rela} and Proposition \ref{propmccmQ}, we obtain that
$\Psi$ is injective.\ Consider a $n$-oscillating tableau $Q$ of length
$l$.\ Set $x_{1}=1$ and for $i=2,...,l$, $x_{i}=k$ if $Q_{i}$ differs from
$Q_{i-1}$ by adding a box in row $k$ and $x_{i}=\overline{k}$ if $Q_{i}$
differs from $Q_{i-1}$ by removing a box in row $k.\;$Consider $w_{Q}%
=x_{l}\cdot\cdot\cdot x_{2}1$.\ Then $Q(w_{Q})=Q$.\ By Theorem \ref{thKN}, the
image of $B(w_{Q})$ by $\Psi$ consists in the pairs $(P,Q)$ where $P$ is a
symplectic tableau of shape $Q_{l}$.\ We deduce immediately that $\Psi$ is surjective.
\end{proof}

In our correspondence, the inverse bumping algorithm is implicit. It is
possible to describe it explicitly but it would be rather cumbersome to do.

\section{A sliding algorithm for $U_{q}(sp_{2n})\label{Sec_JDT}$}

This section is concerned with a symplectic Jeu de Taquin (or sliding
algorithm) introduced by J. T. Sheats \cite{SH} in order to obtain an explicit
bijection between King's and De Concini's symplectic tableaux. The sliding
algorithm for $U_{q}(sl_{n})$ is confluent because its steps stay within the
initial plactic congruence class (see \cite{Fu}). Sheats has conjectured this
property for his Jeu de Taquin. Our aim is now to prove this conjecture and
thus obtain an alternative way to compute $P(w)$ for any word $w.$

We first recall the ideas of Sheats and translate them into
Kashiwara-Nakashima's conventions. Next we extend his algorithm in order to
make it compatible with the relation of contraction $R_{3}.$ Finally we show
that this algorithm is also confluent.

\subsection{Sheats sliding algorithm}

\subsubsection{Skew admissible tableaux}

Let $\lambda=\overset{n}{\underset{i=1}{\sum}}\lambda_{i}\Lambda_{i}$ and
$\mu=\overset{n}{\underset{i=1}{\sum}}\mu_{i}\Lambda_{i}$ be two dominant
weights such that $\mu_{i}\leq\lambda_{i}$ for $i=1,...,n$.\ A skew tableau of
shape $\lambda/\mu$ over $\mathcal{C}_{n}$ is a filling of the skew Young
diagram $Y_{\lambda}/Y_{\mu}$ by letters of $\mathcal{C}_{n}$ giving columns
strictly decreasing from top to bottom.

\begin{definition}
\label{def_skew_tab}A skew tableau over $\mathcal{C}_{n}$ is admissible if its
columns are admissible and the rows of its split form (obtained by turning
each column $C$ into its split form $(lC,rC)$) are weakly increasing from left
to right.
\end{definition}

\begin{example}
$T=%
\begin{tabular}
[c]{ll|ll}\cline{3-4}%
&  & $\mathtt{4}$ & \multicolumn{1}{|l|}{$\mathtt{4}$}\\\cline{2-4}%
& \multicolumn{1}{|l|}{$\mathtt{3}$} & $\mathtt{\bar{4}}$ &
\multicolumn{1}{|l|}{$\mathtt{\bar{3}}$}\\\hline
\multicolumn{1}{|l}{$\mathtt{\bar{3}}$} & \multicolumn{1}{|l|}{$\mathtt{\bar
{3}}$} & $\mathtt{\bar{2}}$ & \multicolumn{1}{|l}{}\\\cline{1-3}%
\multicolumn{1}{|l}{$\mathtt{\bar{2}}$} & \multicolumn{1}{|l|}{$\mathtt{\bar
{2}}$} &  & \\\cline{1-2}%
\end{tabular}
$ is a admissible skew tableau. Its split form is
\[%
\begin{tabular}
[c]{llll|llll}\cline{5-8}%
&  &  &  & $\mathtt{4}$ & \multicolumn{1}{|l}{$\mathtt{4}$} &
\multicolumn{1}{|l}{$\mathtt{4}$} & \multicolumn{1}{|l|}{$\mathtt{4}$%
}\\\cline{3-8}%
&  & \multicolumn{1}{|l}{$\mathtt{1}$} & \multicolumn{1}{|l|}{$\mathtt{3}$} &
$\mathtt{\bar{4}}$ & \multicolumn{1}{|l}{$\mathtt{\bar{3}}$} &
\multicolumn{1}{|l}{$\mathtt{\bar{3}}$} & \multicolumn{1}{|l|}{$\mathtt{\bar
{3}}$}\\\hline
\multicolumn{1}{|l}{$\mathtt{\bar{3}}$} & \multicolumn{1}{|l}{$\mathtt{\bar
{3}}$} & \multicolumn{1}{|l}{$\mathtt{\bar{3}}$} &
\multicolumn{1}{|l|}{$\mathtt{\bar{2}}$} & $\mathtt{\bar{2}}$ &
\multicolumn{1}{|l}{$\mathtt{\bar{2}}$} & \multicolumn{1}{|l}{} &
\\\cline{1-4}\cline{1-6}%
\multicolumn{1}{|l}{$\mathtt{\bar{2}}$} & \multicolumn{1}{|l}{$\mathtt{\bar
{2}}$} & \multicolumn{1}{|l}{$\mathtt{\bar{2}}$} &
\multicolumn{1}{|l|}{$\mathtt{\bar{1}}$} &  &  &  & \\\cline{1-4}%
\end{tabular}
\]
\end{example}

\begin{lemma}
\label{Lem_Sktab_stable_by_op}The set of readings of the admissible skew
tableaux of shape $\lambda/\mu$ is a sub-crystal of $G_{n}$.
\end{lemma}

\begin{proof}
Let $r<n$. Denote by $U_{r}$ the subalgebra of $U_{q}(sp_{2n})$ generated by
$\widetilde{e}_{i},\ \widetilde{f}_{i}$, $i=r+1,...,n$.\ Clearly, $U_{r}\cong
U_{q}(sp_{2(n-r)})$.\ By restriction to $U_{r},$ any $U_{q}(sp_{2n})$-crystal
gives a $U_{r}$-crystal obtained by forgetting the arrows of color
$i=1,...,r$. In particular, $B(\lambda)$ decomposes into a union of connected
$U_{r}$-crystals whose vertices are labelled by symplectic tableaux of shape
$\lambda$.\ Moreover, for all tableaux of such a connected component, the
subtableau consisting of the letters $1,2,...,r$ is the same since
$\widetilde{e}_{i},\ \widetilde{f}_{i}$ $(i>r)$ leaves these letters
unchanged. Hence if we consider the set $B(\lambda,\mu)\subset B(\lambda)$ of
readings of the symplectic tableaux $T$ of shape $\lambda$ on $\mathcal{C}%
_{n}$ such that:

\noindent$\mathrm{(i)}$ the sub tableau of $T$ formed by the letters
$i=1,...,r$ is the fixed tableau of shape $\mu$ having, for each $i,$ only
letters $i$ on its $i$-th row, and

\noindent$\mathrm{(ii)}$ $T$ does not contain any letters in $\{\overline
{r},...,\overline{1}\},$

\noindent then we see that $B(\lambda,\mu)$ is stable under $\widetilde{e}%
_{i},\ \widetilde{f}_{i}$ $(i>r),$ hence is a union of $U_{r}$%
-crystals.\ Finally, this set can be identified to the set of readings of all
admissible skew tableaux of shape $\lambda/\mu$ over $\{r+1,...,n,\overline
{n},...,\overline{r+1}\}$.\ Hence shifting the indices by $r$, we have proved
that the set of readings of the admissible skew tableaux of shape $\lambda
/\mu$ over $\{1,...,n-r,\overline{n-r},...,\overline{1}\}$ is a sub crystal of
$G_{n-r}$ and since $r<n$ are arbitrary, we are done.
\end{proof}

This lemma implies that the readings of the admissible skew tableaux of shape
$\lambda/\mu$ are the vertices of a sub-crystal of $G_{n}$. We denote by
$\mathcal{T}_{(\lambda/\mu)}$ the set of admissible skew tableaux of shape
$(\lambda/\mu)$ and by $\mathcal{U}_{(\lambda/\mu)}$ the set of readings of
these skew tableaux.

\begin{definition}
Consider an admissible skew tableau of shape $\lambda/\mu$. An inner corner is
a box of $Y_{\mu}$ such that the boxes down and to the right are not in
$Y_{\mu}.$ An outside corner is a box of $Y_{\lambda}$ such that the boxes
down and to the right are not in $Y_{\lambda}.$
\end{definition}

\begin{definition}
A skew tableau is punctured if one of its box contains the symbol $\ast$
called the puncture.

A punctured column $C$ is admissible if the column $C^{\prime}$ obtained by
ignoring the puncture is admissible.\ Then the punctured columns $rC$ and $lC$
are respectively obtained by replacing the letters of $C$ (except the
puncture) by the letters of $rC^{\prime}$ and $lC^{\prime}$.\ The split form
of $C$ is $lCrC.$

A punctured skew tableau is admissible if its columns are admissible and the
rows of its split form (obtained by splitting its columns) are weakly
increasing (ignoring the puncture).
\end{definition}

\begin{example}
$T=%
\begin{tabular}
[c]{ll|ll}\cline{3-4}%
&  & $\mathtt{4}$ & \multicolumn{1}{|l|}{$\mathtt{4}$}\\\cline{2-4}%
& \multicolumn{1}{|l|}{$\mathtt{3}$} & $\mathtt{\ast}$ &
\multicolumn{1}{|l|}{$\mathtt{\bar{3}}$}\\\hline
\multicolumn{1}{|l}{$\mathtt{\bar{3}}$} & \multicolumn{1}{|l|}{$\mathtt{\bar
{3}}$} & $\mathtt{\bar{2}}$ & \multicolumn{1}{|l}{}\\\cline{1-3}%
\multicolumn{1}{|l}{$\mathtt{\bar{2}}$} & \multicolumn{1}{|l|}{$\mathtt{\bar
{2}}$} &  & \\\cline{1-2}%
\end{tabular}
$ is a admissible skew punctured tableau of split form $\vspace{0.2cm}spl(T)=%
\begin{tabular}
[c]{llll|llll}\cline{5-8}%
&  &  &  & $\mathtt{3}$ & \multicolumn{1}{|l}{$\mathtt{4}$} &
\multicolumn{1}{|l}{$\mathtt{4}$} & \multicolumn{1}{|l|}{$\mathtt{4}$%
}\\\cline{3-8}%
&  & \multicolumn{1}{|l}{$\mathtt{1}$} & \multicolumn{1}{|l|}{$\mathtt{3}$} &
$\mathtt{\ast}$ & \multicolumn{1}{|l}{$\mathtt{\ast}$} &
\multicolumn{1}{|l}{$\mathtt{\bar{3}}$} & \multicolumn{1}{|l|}{$\mathtt{\bar
{3}}$}\\\hline
\multicolumn{1}{|l}{$\mathtt{\bar{3}}$} & \multicolumn{1}{|l}{$\mathtt{\bar
{3}}$} & \multicolumn{1}{|l}{$\mathtt{\bar{3}}$} &
\multicolumn{1}{|l|}{$\mathtt{\bar{2}}$} & $\mathtt{\bar{2}}$ &
\multicolumn{1}{|l}{$\mathtt{\bar{2}}$} & \multicolumn{1}{|l}{} &
\\\cline{1-4}\cline{1-6}%
\multicolumn{1}{|l}{$\mathtt{\bar{2}}$} & \multicolumn{1}{|l}{$\mathtt{\bar
{2}}$} & \multicolumn{1}{|l}{$\mathtt{\bar{2}}$} &
\multicolumn{1}{|l|}{$\mathtt{\bar{1}}$} &  &  &  & \\\cline{1-4}%
\end{tabular}
$.
\end{example}

\subsubsection{Elementary step of the sliding algorithm
\label{subsub_sec_fund_step}}

Let us consider $T$ an admissible punctured skew tableau containing two
columns $C_{1}$ and $C_{2}$ with the puncture in $C_{1}.$ To apply an
elementary step of the sliding algorithm to $T$ we have first to consider the
split form of $T$. In this split form we have a configuration of the type:%

\begin{tabular}
[c]{|l|l|ll}\cline{1-2}%
$...$ & $...$ & $...$ & $...$\\\hline
$\ast$ & $\ast$ & $b$ & \multicolumn{1}{|l|}{$b^{\prime}$}\\\hline
$a$ & $a^{\prime}$ & $...$ & \multicolumn{1}{|l|}{$...$}\\\hline
$...$ & $...$ &  & \multicolumn{1}{|l|}{}\\\hline
\end{tabular}
where the boxes containing $a,a^{\prime}$ and $b,b^{\prime}$ may be
empty.\newline An elementary step of the Symplectic Jeu de Taquin
($\mathrm{SJDT}$) consists of the following transformations:

\begin{enumerate}
\item  If $a^{\prime}\leq b$ or the double box $b$ $b^{\prime}$ is empty, then
the doubles boxes $a$ $a^{\prime}$ and $\ast$ $\ast$ are permuted

\item  If $a^{\prime}>b$ or the double box $a$ $a^{\prime}$ is empty then:

\textrm{(i)}: when $b$ is a barred letter, $b$ slides into $rC_{1}$ to the box
containing $\ast$ and $D_{1}=\Phi(C_{1})-\{\ast\}+\{b\}$ is a co-admissible
column (see section 2.2, (\ref{phi})).\ Simultaneously the symbol $\ast$
slides into $lC_{2}$ to the box containing $b$ and $C_{2}^{\prime}%
=C_{2}-\{b\}+\{\ast\}$ is a punctured admissible column . Then we obtain a new
punctured skew tableau $C_{1}^{\prime}C_{2}^{\prime}$ by setting
$C_{1}^{\prime}=\Phi^{-1}(D_{1})$.

($\mathrm{ii}$): when $b$ is an unbarred letter, $b$ slides into $rC_{1}$ to
the box containing $\ast$ and give a new column $C_{1}^{\prime}=C_{1}%
-\{\ast\}+\{b\}$.\ Simultaneously the symbol $\ast$ slides into $lC_{2}$ to
the box containing $b$ and $D_{2}=\Phi(C_{2})-\{b\}+\{\ast\}$ is a punctured
co-admissible column. Then we obtain a new punctured skew tableau
$C_{1}^{\prime}C_{2}^{\prime}$ by setting $C_{2}^{\prime}=\Phi^{-1}(D_{1})$.
\end{enumerate}

Notice that in case $2$ \textrm{(i) }the co-admissibility of $D_{1}$ is not
immediate and in case $2$ \textrm{(ii) }the column $C_{1}^{\prime}$ may be not admissible.

\begin{lemma}
\label{lemmfunJDTS} We can always apply an elementary step of the
$\mathrm{SJDT}$ to an admissible punctured skew symplectic tableau $($i.e.
$D_{1}$ is a co-admissible column in case $2$ $\mathrm{(i)})$.
\end{lemma}

\begin{proof}
See Lemma 9.1 of \cite{SH}.
\end{proof}

\begin{example}
\ \ \ \ \ \ \ \ \ \label{ex_horizontal_move}

\noindent For $T_{1}=%
\begin{tabular}
[c]{|l|l}\hline
$\mathtt{2}$ & \multicolumn{1}{|l|}{$\mathtt{4}$}\\\hline
$\mathtt{4}$ & \multicolumn{1}{|l|}{$\mathtt{5}$}\\\hline
$\mathtt{{{{{{{{{{{{{{{{{{{\ast}}}}}}}}}}}}}}}}}}}$ &
\multicolumn{1}{|l|}{$\mathtt{\bar{4}}$}\\\hline
$\mathtt{\bar{3}}$ & \multicolumn{1}{|l|}{$\mathtt{\bar{1}}$}\\\hline
$\mathtt{\bar{1}}$ & \\\cline{1-1}%
\end{tabular}
$ $spl(T_{1})=%
\begin{tabular}
[c]{|l|l|ll}\hline
$\mathtt{2}$ & $\mathtt{2}$ & $\mathtt{3}$ & \multicolumn{1}{|l|}{$\mathtt{4}%
$}\\\hline
$\mathtt{4}$ & $\mathtt{4}$ & $\mathtt{5}$ & \multicolumn{1}{|l|}{$\mathtt{5}%
$}\\\hline
$\mathtt{{{{{{{{{{{{{{{{{{{\ast}}}}}}}}}}}}}}}}}}}$ & $\mathtt{\ast}$ &
$\mathtt{\bar{4}}$ & \multicolumn{1}{|l|}{$\mathtt{\bar{3}}$}\\\hline
$\mathtt{\bar{3}}$ & $\mathtt{\bar{3}}$ & $\mathtt{\bar{1}}$ &
\multicolumn{1}{|l|}{$\mathtt{\bar{1}}$}\\\hline
$\mathtt{\bar{1}}$ & $\mathtt{\bar{1}}$ &  & \\\cline{1-2}%
\end{tabular}
$. We are in case $2$ $(\mathrm{i})$ and $C_{1}^{\prime}C_{2}^{\prime}=$
\begin{tabular}
[c]{|l|l}\hline
$\mathtt{2}$ & \multicolumn{1}{|l|}{$\mathtt{4}$}\\\hline
$\mathtt{5}$ & \multicolumn{1}{|l|}{$\mathtt{5}$}\\\hline
$\mathtt{\bar{5}}$ & \multicolumn{1}{|l|}{$\mathtt{\ast}$}\\\hline
$\mathtt{\bar{3}}$ & \multicolumn{1}{|l|}{$\mathtt{\bar{1}}$}\\\hline
$\mathtt{\bar{1}}$ & \\\cline{1-1}%
\end{tabular}
.

\noindent For $T_{2}=%
\begin{tabular}
[c]{|l|l}\hline
$\mathtt{2}$ & \multicolumn{1}{|l|}{$\mathtt{2}$}\\\hline
$\mathtt{3}$ & \multicolumn{1}{|l|}{$\mathtt{3}$}\\\hline
$\mathtt{{{{{{{{{{{{{{{{{{{\ast}}}}}}}}}}}}}}}}}}}$ &
\multicolumn{1}{|l|}{$\mathtt{5}$}\\\hline
$\mathtt{\bar{5}}$ & \multicolumn{1}{|l|}{$\mathtt{\bar{5}}$}\\\hline
$\mathtt{\bar{1}}$ & \\\cline{1-1}%
\end{tabular}
$ $spl(T_{2})=%
\begin{tabular}
[c]{|l|l|ll}\hline
$\mathtt{2}$ & $\mathtt{2}$ & $\mathtt{2}$ & \multicolumn{1}{|l|}{$\mathtt{2}%
$}\\\hline
$\mathtt{3}$ & $\mathtt{3}$ & $\mathtt{3}$ & \multicolumn{1}{|l|}{$\mathtt{3}%
$}\\\hline
$\mathtt{{{{{{{{{{{{{{{{{{{\ast}}}}}}}}}}}}}}}}}}}$ & $\mathtt{\ast}$ &
$\mathtt{4}$ & \multicolumn{1}{|l|}{$\mathtt{5}$}\\\hline
$\mathtt{\bar{5}}$ & $\mathtt{\bar{5}}$ & $\mathtt{\bar{5}}$ &
\multicolumn{1}{|l|}{$\mathtt{\bar{4}}$}\\\hline
$\mathtt{\bar{1}}$ & $\mathtt{\bar{1}}$ &  & \\\cline{1-2}%
\end{tabular}
$. We are in case $2$ $(\mathrm{ii})$ and $C_{1}^{\prime}C_{2}^{\prime}=$
\begin{tabular}
[c]{|l|l}\hline
$\mathtt{2}$ & \multicolumn{1}{|l|}{$\mathtt{2}$}\\\hline
$\mathtt{3}$ & \multicolumn{1}{|l|}{$\mathtt{3}$}\\\hline
$\mathtt{4}$ & \multicolumn{1}{|l|}{$\mathtt{\ast}$}\\\hline
$\mathtt{\bar{5}}$ & \multicolumn{1}{|l|}{$\mathtt{\bar{4}}$}\\\hline
$\mathtt{\bar{1}}$ & \\\cline{1-1}%
\end{tabular}
.\vspace{0.2cm}

\noindent For $T_{3}=%
\begin{tabular}
[c]{|l|l|l|}\hline
$\mathtt{4}$ & $\mathtt{\ast}$ & $\mathtt{4}$\\\hline
$\mathtt{\bar{5}}$ & $\mathtt{\bar{4}}$ & $\mathtt{\bar{3}}$\\\hline
\end{tabular}
$ we obtain $%
\begin{tabular}
[c]{|l|l|l|}\hline
$\mathtt{4}$ & $\mathtt{4}$ & $\mathtt{\ast}$\\\hline
$\mathtt{\bar{5}}$ & $\mathtt{\bar{4}}$ & $\mathtt{\bar{3}}$\\\hline
\end{tabular}
$.
\end{example}

The punctured skew tableau obtained by computing a step of the $\mathrm{SJDT}$
on an admissible skew tableau is not always admissible. In the second example
above $C_{1}^{\prime}$ is not admissible and in the third the rows of the
split form are not increasing (we will see that this last problem does not
occur in the complete $\mathrm{SJDT}$ algorithm).

\subsubsection{Complete symplectic Jeu de Taquin ($\mathrm{SJDT}$)}

Let $T$ be an admissible skew tableau and $c$ an inner corner in $T.$ In order
to apply the complete sliding algorithm let us puncture the corner $c.$ We
obtain an admissible punctured skew tableau. To see what happens when we apply
successively elementary steps of $\mathrm{SJDT}$ to this skew tableau, we need
to compute the split form for each intermediate punctured tableau. We have
seen that a horizontal move of an unbarred letter may give a new non
admissible column $C_{1}^{\prime}$ such that all the strict subwords of
$\mathrm{w}(C_{1}^{\prime})$ are admissible.\ So it is impossible to compute
its split form using letters of $\mathcal{C}_{n}$. To overcome this problem,
we embed the alphabet $\mathcal{C}_{n}$ into%
\[
\mathcal{C}_{n+1}^{\prime}=\{a_{1}<1<\cdot\cdot\cdot<n<\overline{n}<\cdot
\cdot\cdot<\overline{1}<\overline{a_{1}}\}.
\]
To compute the split form of a non admissible column $C$ such that all the
strict subwords of $\mathrm{w}(C)$ are admissible, we extend the algorithm of
Subsection \ref{subsec_splitC} by using the new letter $a_{1}.$ For example%
\[
\text{if }C=%
\begin{tabular}
[c]{|l|}\hline
\texttt{2}\\\hline
\texttt{3}\\\hline
\texttt{4}\\\hline
\texttt{\={4}}\\\hline
\texttt{\={1}}\\\hline
\end{tabular}
\text{ }(lC,rC)=%
\begin{tabular}
[c]{|l|l|}\hline
$a_{1}$ & \texttt{2}\\\hline
\texttt{2} & \texttt{3}\\\hline
\texttt{3} & \texttt{4}\\\hline
\texttt{\={4}} & \texttt{\={1}}\\\hline
\texttt{\={1}} & $\overline{a_{1}}$\\\hline
\end{tabular}
\text{ in }\mathcal{C}_{n+1}^{\prime}\text{.}%
\]
So all the columns that may be obtained when we apply an elementary step of
$\mathrm{SJDT}$ to an admissible skew tableau (defined on $\mathcal{C}_{n}$)
can be split in $\mathcal{C}_{n+1}^{\prime}$.\ We will say that a skew
punctured tableau is $a_{1}$-admissible if all its columns can be split in
$\mathcal{C}_{n+1}^{\prime}$ and the rows of the obtained split form are
weakly increasing.

\begin{theorem}
(Sheats \cite{SH})\label{THJDT}

\begin{itemize}
\item  Elementary steps of $\mathrm{SJDT}$ may be applied to $T$ until the
puncture $\ast$ becomes an outside corner.

\item  All the skew punctured tableaux obtained as steps in the algorithm are
$a_{1}$-admissible. Moreover $\overline{a_{1}}$ and $a_{1}$ may only appear
simultaneously in the split form of the column containing the inner corner $c$
of $T$ at which the slide started.
\end{itemize}
\end{theorem}

\begin{proof}
See Proposition 9.2 of \cite{SH}.
\end{proof}

\begin{example}
Suppose $spl(T)=$
\begin{tabular}
[c]{ll|l|l|ll}\cline{3-6}%
&  & $\mathtt{\ast}$ & $\mathtt{\ast}$ & $\mathtt{2}$ &
\multicolumn{1}{|l|}{$\mathtt{2}$}\\\cline{3-6}%
&  & $\mathtt{2}$ & $\mathtt{2}$ & $\mathtt{3}$ &
\multicolumn{1}{|l|}{$\mathtt{3}$}\\\hline
\multicolumn{1}{|l}{$\mathtt{2}$} & \multicolumn{1}{|l|}{$\mathtt{3}$} &
$\mathtt{3}$ & $\mathtt{3}$ & $\mathtt{4}$ & \multicolumn{1}{|l|}{$\mathtt{4}%
$}\\\hline
\multicolumn{1}{|l}{$\mathtt{\bar{5}}$} & \multicolumn{1}{|l|}{$\mathtt{\bar
{5}}$} & $\mathtt{\bar{4}}$ & $\mathtt{\bar{4}}$ & $\mathtt{\bar{1}}$ &
\multicolumn{1}{|l|}{$\mathtt{\bar{1}}$}\\\hline
\multicolumn{1}{|l}{$\mathtt{\bar{3}}$} & \multicolumn{1}{|l|}{$\mathtt{\bar
{2}}$} & $\mathtt{\bar{1}}$ & $\mathtt{\bar{1}}$ &  & \\\cline{1-4}%
\end{tabular}
. We compute successively\vspace{0.5cm} the split form of the $a_{1}%
$-admissible punctured skew tableaux:\vspace{0.5cm}

$%
\begin{tabular}
[c]{ll|l|l|ll}\cline{3-6}%
&  & $\mathtt{2}$ & $\mathtt{2}$ & $\mathtt{2}$ &
\multicolumn{1}{|l|}{$\mathtt{2}$}\\\cline{3-6}%
&  & $\mathtt{\ast}$ & $\mathtt{\ast}$ & $\mathtt{3}$ &
\multicolumn{1}{|l|}{$\mathtt{3}$}\\\hline
\multicolumn{1}{|l}{$\mathtt{2}$} & \multicolumn{1}{|l|}{$\mathtt{3}$} &
$\mathtt{3}$ & $\mathtt{3}$ & $\mathtt{4}$ & \multicolumn{1}{|l|}{$\mathtt{4}%
$}\\\hline
\multicolumn{1}{|l}{$\mathtt{\bar{5}}$} & \multicolumn{1}{|l|}{$\mathtt{\bar
{5}}$} & $\mathtt{\bar{4}}$ & $\mathtt{\bar{4}}$ & $\mathtt{\bar{1}}$ &
\multicolumn{1}{|l|}{$\mathtt{\bar{1}}$}\\\hline
\multicolumn{1}{|l}{$\mathtt{\bar{3}}$} & \multicolumn{1}{|l|}{$\mathtt{\bar
{2}}$} & $\mathtt{\bar{1}}$ & $\mathtt{\bar{1}}$ &  & \\\cline{1-4}%
\end{tabular}
\,;$\hspace{1.3cm}$\mathtt{%
\begin{tabular}
[c]{ll|l|l|ll}\cline{3-6}%
&  & $\mathtt{2}$ & $\mathtt{2}$ & $\mathtt{2}$ &
\multicolumn{1}{|l|}{$\mathtt{2}$}\\\cline{3-6}%
&  & $\mathtt{3}$ & $\mathtt{3}$ & $\mathtt{3}$ &
\multicolumn{1}{|l|}{$\mathtt{3}$}\\\hline
\multicolumn{1}{|l}{$\mathtt{2}$} & \multicolumn{1}{|l|}{$\mathtt{3}$} &
$\mathtt{\ast}$ & $\mathtt{\ast}$ & $\mathtt{4}$ &
\multicolumn{1}{|l|}{$\mathtt{4}$}\\\hline
\multicolumn{1}{|l}{$\mathtt{\bar{5}}$} & \multicolumn{1}{|l|}{$\mathtt{\bar
{5}}$} & $\mathtt{\bar{4}}$ & $\mathtt{\bar{4}}$ & $\mathtt{\bar{1}}$ &
\multicolumn{1}{|l|}{$\mathtt{\bar{1}}$}\\\hline
\multicolumn{1}{|l}{$\mathtt{\bar{3}}$} & \multicolumn{1}{|l|}{$\mathtt{\bar
{2}}$} & $\mathtt{\bar{1}}$ & $\mathtt{\bar{1}}$ &  & \\\cline{1-4}%
\end{tabular}
}\,;\hspace{1.3cm}\mathtt{%
\begin{tabular}
[c]{ll|l|l|ll}\cline{3-6}%
&  & $a_{1}$ & $\mathtt{2}$ & $\mathtt{2}$ & \multicolumn{1}{|l|}{$\mathtt{2}%
$}\\\cline{3-6}%
&  & $\mathtt{2}$ & $\mathtt{3}$ & $\mathtt{3}$ &
\multicolumn{1}{|l|}{$\mathtt{3}$}\\\hline
\multicolumn{1}{|l}{$\mathtt{2}$} & \multicolumn{1}{|l|}{$\mathtt{3}$} &
$\mathtt{3}$ & $\mathtt{4}$ & $\mathtt{\ast}$ &
\multicolumn{1}{|l|}{$\mathtt{\ast}$}\\\hline
\multicolumn{1}{|l}{$\mathtt{\bar{5}}$} & \multicolumn{1}{|l|}{$\mathtt{\bar
{5}}$} & $\mathtt{\bar{4}}$ & $\mathtt{\bar{1}}$ & $\mathtt{\bar{1}}$ &
\multicolumn{1}{|l|}{$\mathtt{\bar{1}}$}\\\hline
\multicolumn{1}{|l}{$\mathtt{\bar{3}}$} & \multicolumn{1}{|l|}{$\mathtt{\bar
{2}}$} & $\mathtt{\bar{1}}$ & $\overline{a_{1}}$ &  & \\\cline{1-4}%
\end{tabular}
}$\texttt{\vspace{0.5cm}}\linebreak \texttt{%
\begin{tabular}
[c]{ll|l|l|ll}\cline{3-6}%
&  & $a_{1}$ & $\mathtt{2}$ & $\mathtt{2}$ & \multicolumn{1}{|l|}{$\mathtt{2}%
$}\\\cline{3-6}%
&  & $\mathtt{2}$ & $\mathtt{3}$ & $\mathtt{3}$ & \multicolumn{1}{|l|}{$3$%
}\\\hline
\multicolumn{1}{|l}{$\mathtt{2}$} & \multicolumn{1}{|l|}{$\mathtt{3}$} &
$\mathtt{3}$ & $\mathtt{4}$ & $\mathtt{\bar{1}}$ & \multicolumn{1}{|l|}{$\bar
{1}$}\\\hline
\multicolumn{1}{|l}{$\mathtt{\bar{5}}$} & \multicolumn{1}{|l|}{$\mathtt{\bar
{5}}$} & $\mathtt{\bar{4}}$ & $\mathtt{\bar{1}}$ & $\mathtt{\ast}$ &
\multicolumn{1}{|l|}{$\mathtt{\ast}$}\\\hline
\multicolumn{1}{|l}{$\mathtt{\bar{3}}$} & \multicolumn{1}{|l|}{$\mathtt{\bar
{2}}$} & $\mathtt{\bar{1}}$ & $\overline{a_{1}}$ &  & \\\cline{1-4}%
\end{tabular}
}. Then we obtain the $a_{1}$-admissible skew tableau: \texttt{%
\begin{tabular}
[c]{l|l|l}\cline{2-3}%
& $\mathtt{2}$ & \multicolumn{1}{|l|}{$\mathtt{2}$}\\\cline{2-3}%
& $\mathtt{3}$ & \multicolumn{1}{|l|}{$\mathtt{3}$}\\\hline
\multicolumn{1}{|l|}{$\mathtt{3}$} & $\mathtt{4}$ &
\multicolumn{1}{|l|}{$\mathtt{\bar{1}}$}\\\hline
\multicolumn{1}{|l|}{$\mathtt{\bar{5}}$} & $\mathtt{\bar{4}}$ & \\\cline{1-2}%
\multicolumn{1}{|l|}{$\mathtt{\bar{3}}$} & $\mathtt{\bar{1}}$ & \\\cline{1-2}%
\end{tabular}
}.
\end{example}

\subsection{Sliding algorithm on $\mathcal{C}_{n}$}

Let $T$ be an admissible skew tableau and $c$ an inner corner. If we denote by
$T^{\prime}$ the skew tableau obtained by applying the complete $\mathrm{SJDT}%
$ to $T$, then $T^{\prime}$ may be only $a_{1}$-admissible (see Theorem
\ref{THJDT}). Suppose that, in the split form, $\overline{a_{1}}$ and $a_{1}$
occur in the $k$-th split column $lC_{k}^{\prime}\,rC_{k}^{\prime}$ of
$T^{\prime}.$ Then the column $C_{k}^{\prime}$ is not admissible and, in order
to obtain an admissible skew tableau, we are led to consider the skew tableau
$\widetilde{T}^{\prime}$ obtained by erasing the top and bottom boxes of
$C_{k}^{\prime}$ and filling this new column with the letters of the word
$\widetilde{\mathrm{w(}C_{k}^{\prime})}$. We denote this new column by
$\widetilde{C}_{k}$.

\begin{example}
Continuing the previous example we obtain:%
\[
\widetilde{T}=%
\begin{tabular}
[c]{lll}\cline{3-3}%
&  & \multicolumn{1}{|l|}{$\mathtt{2}$}\\\cline{2-3}%
& \multicolumn{1}{|l}{$\mathtt{2}$} & \multicolumn{1}{|l|}{$\mathtt{3}$%
}\\\hline
\multicolumn{1}{|l}{$\mathtt{3}$} & \multicolumn{1}{|l}{$\mathtt{3}$} &
\multicolumn{1}{|l|}{$\mathtt{\bar{1}}$}\\\hline
\multicolumn{1}{|l}{$\mathtt{\bar{5}}$} & \multicolumn{1}{|l}{$\mathtt{\bar
{1}}$} & \multicolumn{1}{|l}{}\\\cline{1-2}%
\multicolumn{1}{|l}{$\mathtt{\bar{2}}$} & \multicolumn{1}{|l}{} &
\\\cline{1-1}%
\end{tabular}
.
\]
\end{example}

Using the notations introduced above, we have:

\begin{proposition}
\label{propcontJDTS}$\widetilde{T}^{\prime}$ is an admissible skew tableau and
$\mathrm{w(}T^{\prime})\equiv\mathrm{w(}\widetilde{T}^{\prime})$.
\end{proposition}

\begin{proof}
By Theorem \ref{THJDT}, we know that $T^{\prime}$ is $a_{1}$-admissible and
$(a_{1},\overline{a_{1}})$ occurs only in $C_{k}^{\prime}$. This implies that
there is no box of $T^{\prime}$ weakly below and strictly to the right of the
box containing $\overline{a_{1}}$. So $\widetilde{T}$ has the shape of a skew
tableau. The pair $(a_{1},\overline{a_{1}})$ disappears when we compute
$(l\widetilde{C}_{k},r\widetilde{C}_{k})$. Moreover each letter $x\neq
\overline{a_{1}}$ of $rC_{k}^{\prime}$ is turned into a letter $y\leq x$ in
$r\widetilde{C}_{k}$ and each letter $t\neq a_{1}$ of $lC_{k}^{\prime}$ is
turned into a letter $v\geq t$ in $l\widetilde{C}_{k}$. Thus the rows of the
split form of $\widetilde{T}^{\prime}$ are weakly increasing: $\widetilde
{T}^{\prime}$ is an admissible skew tableau. Finally, the identity
\textrm{w(}$T^{\prime})\equiv_{n}$\textrm{w(}$\widetilde{T}^{\prime})$ is
clear because \textrm{w(}$C_{k}^{\prime})\equiv_{n}$\textrm{w(}$\widetilde
{C}_{k})$ by $R_{3}$.
\end{proof}

Given an admissible skew tableau $T$ and $c$ an inside corner in $T$ we can
apply elementary steps of $\mathrm{SJDT}$ to obtain a skew tableau $T^{\prime
}$.\ We set:%
\[
\mathrm{SJDT}(T,c)=\left\{
\begin{tabular}
[c]{l}%
$T^{\prime}$ if $T^{\prime}$ is admissible\\
$\widetilde{T}$ if $T^{\prime}$ is only $a$-admissible
\end{tabular}
\right.  .
\]
During the algorithm an inner corner is filled or $\mathrm{SJDT}(c,T)$ has two
boxes less than $T$. By choosing a new inner corner at each step, we can
iterate the procedure $T\rightarrow\mathrm{SJDT}(T,c)$ to construct a
symplectic tableau from any admissible skew tableau. Now we are going to
establish the

\begin{conjecture}
(Sheats) The symplectic tableau obtained by iterating the sliding algorithm is
independent of the order in which the inner corners are filled.
\end{conjecture}

\subsection{Proof of Sheats' conjecture}

In this section we show that, for any fixed inner corner $c$, the map
$\mathrm{w}(T)\rightarrow\mathrm{w}(\mathrm{SJDT}(T,c))$ defined on
$U_{(\lambda/\mu)}$ is a crystal isomorphism. This result will imply Sheats' conjecture.

We begin with skew tableaux of two columns. Consider $p,q,k$ some integers
such that $0<k\leq p\leq q$ and denote by $\mathcal{T}_{(q,p)/k}$ the set of
admissible skew tableaux of two columns of shape $\lambda/\mu$ where
$Y_{\lambda}$ consists in the two columns of height $q$ and $p$ and $Y_{\mu}$
is the single column of height $k$. The skew tableaux of $\mathcal{T}%
_{(q,p)/k}$ have a unique inner corner $c$.\ By Lemma
\ref{Lem_Sktab_stable_by_op}, the subset $\mathcal{U}_{(q,p)/k}$ consisting of
the readings of the skew tableaux of $\mathcal{T}_{(q,p)/k}$ is a sub-crystal
of $G_{n}$.\ Set $\Psi_{(q,p)/k}(0)=0$ and for any $T\in\mathcal{T}_{(q,p)/k}$%
\[
\Psi_{(q,p)/k}(\mathrm{w}(T))=\mathrm{w(SJDT}(T,c)).
\]
Then $\Psi_{(q,p)/k}$ is a map defined on $\mathcal{U}_{(q,p)/k}\cup
\{0\}$.\ We are going to prove that $\Psi_{(q,p)/k}$ is in fact a crystal
isomorphism from $\mathcal{U}_{(q,p)/k}$ to its image.\ To obtain this result
it suffices to see that $\Psi_{(q,p)/k}$ verifies the assertions
$\mathrm{(ii)}$ and $\mathrm{(iii)}$ of Lemma \ref{lem_iso_G}.\ This is the
goal of technical Lemmas \ref{lem_JT_HWV}, \ref{lem_com_en}, \ref{lem_psi+}
and \ref{lem_theta_comm_with_ei}. Elementary steps of $\mathrm{SJDT}$ are
applied on two consecutive columns of skew punctured tableaux. So by using
crystal isomorphisms of type $\Psi_{(q,p)/k}$ and Lemma \ref{lem_util}, it
will be easy to obtain a crystal isomorphism defined on $U_{(\lambda/\mu)}$ in
Theorem \ref{th_iso_JDT}.

We will say that a horizontal move occurs in $T\in\mathcal{T}_{(q,p)/k}$ when
\textrm{w}$(T)\notin\mathcal{U}_{(q-1,p)/(k-1)}$. It means that we have in
$spl(T)$ a configuration
\begin{equation}
H_{(a,b)}=%
\begin{tabular}
[c]{|l|l|l|l|}\hline
$.$ & $.$ & $b$ & $.$\\\hline
$.$ & $a$ & $.$ & $.$\\\hline
\end{tabular}
\label{conf_slide}%
\end{equation}
with $a>b.\;$Notice that $\Psi_{(q,p)/k}(\mathrm{w}(T))=\mathrm{w}(T)$ if no
horizontal move occurs in $T$.

\begin{lemma}
\label{lem_JT_HWV}Let $\mathcal{U}_{(q,p)/k}^{0}$ be the set of highest weight
vertices of $\mathcal{U}_{(q,p)/k}.$ Then $\Psi_{(q,p)/k}$ induces a bijection
from $\mathcal{U}_{(q,p)/k}^{0}$ to $\Psi_{(q,p)/k}(\mathcal{U}_{(q,p)/k}%
^{0})$. Moreover, $\Psi_{(q,p)/k}(\mathcal{U}_{(q,p)/k}^{0})$ contains only
highest weight vertices and for any $w^{0}\in\mathcal{U}_{(q,p)/k}^{0},$
$d(w^{0})=d(\Psi_{(q,p)/k}(w^{0}))$.
\end{lemma}

\begin{proof}
The assertion $d(w^{0})=d(\Psi_{(q,p)/k}(w^{0}))$ for any $w^{0}\in$
$\mathcal{U}_{(q,p)/k}^{0}$ follows immediately from the definition of
$\mathrm{SJDT}$.\ Consider $w_{1},w_{2}\in\mathcal{U}_{(q,p)/k}^{0}$ such that
$\Psi_{(q,p)/k}(w_{1})=\Psi_{(q,p)/k}(w_{2})$.\ Denote respectively by
$C_{1}C_{2}$ and $D_{1}D_{2}$ the skew tableaux of $\mathcal{T}_{(q,p)/k}$ of
reading $w_{1}$ and $w_{2}$.\ If no horizontal move occurs in $C_{1}C_{2}$, no
horizontal move occurs in $D_{1}D_{2}$ and $w_{1}=\Psi_{(q,p)/k}(w_{1}%
)=\Psi_{(q,p)/k}(w_{2})=w_{2}$. So we can suppose that a horizontal move
occurs in $C_{1}C_{2}$ and $D_{1}D_{2}$.\ By Lemma \ref{lmmaximalweight},
$\mathrm{w}(C_{2})$ and $\mathrm{w}(D_{2})$ are vertices of highest weight$.$
Hence $D_{2}=C_{2}$ since $h(C_{2})=h(D_{2})=p$. Write $\mathrm{w}%
(C_{2})=1\cdot\cdot\cdot p$.\ Denote by $b$ the letter of $lC_{2}$ which
slides into $rC_{1}$.\ Suppose that $b<p$.$\;$Then $C_{1}C_{2}$ contains a
configuration%
\[
H_{(a,b)}=%
\begin{tabular}
[c]{|l|l|l|l|}\hline
$.$ & $.$ & $b$ & $b$\\\hline
$.$ & $a$ & $b+1$ & $b+1$\\\hline
\end{tabular}
\]
with $a>b$.\ We have $a\leq b+1$. This implies that $a=b+1$. Then
$\mathrm{w}(C_{1})$ contains the letter $b+1\leq p$. Hence $\mathrm{w}(C_{1})$
contains all the letters $b+1,b,...,1$ because by Lemma \ref{lmmaximalweight},
we must have $\varepsilon_{i}\mathrm{w}(C_{1})=0$ for $i=b,...,1.$ Moreover
$C_{1}$ is admissible so $\overline{b}\notin C_{1}$ hence $b\in rC_{1}$. This
is incompatible with the sliding of $b$ in $rC_{1}$. That means that $p$ is
the letter of $C_{2}$ which slide in $rC_{1}$ when we compute $\Psi
_{(q,p)/k}(w_{1})$.\ Similarly $p$ slides into $rD_{1}$ when we compute
$\Psi_{(q,p)/k}(w_{2})$.\ Set $C_{1}^{\prime}=C_{1}+\{p\}$ and $D_{1}^{\prime
}=D_{1}+\{p\}$.\ If these two columns are admissible we will have
$\Psi_{(q,p)/k}(w_{1})=\mathrm{w}(C_{1}^{\prime}C_{2}^{\prime})$ and
$\Psi_{(q,p)/k}(w_{2})=\mathrm{w}(D_{1}^{\prime}C_{2}^{\prime})$ with
$C_{2}^{\prime}=C_{2}-\{p\}.\;$Indeed the horizontal moves considered are of
type $2$ $\mathrm{(ii)}$. Using $\Psi_{(q,p)/k}(w_{1})=\Psi_{(q,p)/k}(w_{2})$,
we obtain $D_{1}^{\prime}=C_{1}^{\prime}$ hence $D_{1}=C_{1}$. If
$C_{1}^{\prime}$ is not admissible, $C_{2}^{\prime}$ is not admissible for
$\Psi_{(q,p)/k}(w_{1})=\Psi_{(q,p)/k}(w_{2})$.\ Then we have $\widetilde
{\mathrm{w}(C_{1}^{\prime})}=\widetilde{\mathrm{w}(D_{1}^{\prime})}$.\ Denote
respectively by $(z,\overline{z})$ and $(t,\overline{t})$ the pairs of letters
erased in $\mathrm{w}(C_{1}^{\prime})$ and $\mathrm{w}(D_{1}^{\prime})$. Then
$N(z)=z+1$ in $C_{1}^{\prime}$ and $N(t)=t+1$ in $D_{1}^{\prime}$. But
$C_{1}^{\prime}$ and $D_{1}^{\prime}$ may differ only by the pairs
$(z,\overline{z})$ and $(t,\overline{t})$. Hence $N(z)=N(t)$, so $z=t$. We
obtain $C_{1}^{\prime}=D_{1}^{\prime}$ and $C_{1}=D_{1}$.

To prove that $\Psi_{(q,p)/k}(w_{1})$ is of highest weight, it suffices to
show that $\mathrm{w}(C_{1}^{\prime}C_{2}^{\prime})$ is of highest weight
because when $C_{1}^{\prime}$ is not admissible $\Psi_{(q,p)/k}(w_{1}%
)\equiv\mathrm{w}(C_{1}^{\prime}C_{2}^{\prime})$. We have $C_{2}^{\prime
}=1\cdot\cdot\cdot(p-1)$. Hence by Lemma \ref{lmmaximalweight}, $w(C_{1}%
^{\prime}C_{2}^{\prime})$ is of highest weight if and only if
\begin{equation}
\varepsilon_{i}(\mathrm{w}(C_{1}^{\prime}))=0\text{ if }i\neq p-1\text{ and
}\varepsilon_{p-1}(\mathrm{w}(C_{1}^{\prime}))\leq1\text{.}\label{cond2}%
\end{equation}
The columns $C_{1}$ and $C_{1}^{\prime}$ differ only by the letter $p$. Then
(\ref{cond2}) is satisfied for $i\notin\{p-1,p\}$.\ We have $\varepsilon
_{p-1}(\mathrm{w}(C_{1}^{\prime}))\leq1$. For if it were not so,
$\varepsilon_{p-1}(\mathrm{w}(C_{1}^{\prime}))=2$ and $\varepsilon
_{p-1}(\mathrm{w}(C_{1}))=1$ because $C_{1}^{\prime}=C_{1}+\{p\}$. Suppose
$\varepsilon_{p}\mathrm{w}(C_{1}^{\prime})\neq0.$ If $\varepsilon
_{p}(\mathrm{w}(C_{1}))=0$ then $C_{1}$ is of highest weight, we have
$\mathrm{w}(C_{1})=1\cdot\cdot\cdot(q-k)$ and no horizontal move occurs in
$C_{1}C_{2}$. So $\varepsilon_{p}(\mathrm{w}(C_{1}))=1$. Then, with the
notation of (\ref{regle+-}) $\rho_{p}(\mathrm{w}(C_{1}))=-$. So $\rho
(\mathrm{w}(C_{1}^{\prime}))=+-$ because $C_{1}^{\prime}=C_{1}+\{p\}$. Hence
$\varepsilon_{p}(\mathrm{w}(C_{1}^{\prime}))=0$ and we obtain a contradiction.
\end{proof}

\begin{lemma}
\label{lem_com_en}Let $T\in\mathcal{T}_{(q,p)/k}$ and suppose $\widetilde
{f}_{n}\,\mathrm{w}(T)\neq0$. Then $\widetilde{f}_{n}\Psi_{(q,p)/k}%
\,\mathrm{w(}T)=\Psi_{(q,p)/k}\widetilde{f}_{n}\,\mathrm{w}(T).$
\end{lemma}

\begin{proof}
The Lemma is clear if no horizontal move occurs in $T$. In the sequel we
suppose that a horizontal move occurs in $T$.\ Let $T_{n}\in\mathcal{T}%
_{(q,p)/k}$ such that $\mathrm{w(}T_{n})=\widetilde{f}_{n}\,\mathrm{w(}T)$.
Write $T=C_{1}C_{2}$ and $T_{n}=D_{1}D_{2}$. We know that $spl(T)$ contains a
configuration $H_{(a,b)}$ of type (\ref{conf_slide}) such that $b$ slides into
$rC_{1}$ when we apply $\mathrm{SJDT}$ to $T$.\ By considering the\ letters of
$spl(T_{n})$ occurring at the same place as the letters of $H_{(a,b)}$ we
obtain in $spl(T_{n})$ a configuration $H_{(a^{\prime},b^{\prime})}$.\ We have%
\begin{equation}%
\begin{tabular}
[c]{|l|l|l|l|}\hline
$.$ & $.$ & $b$ & $.$\\\hline
$.$ & $a$ & $.$ & $.$\\\hline
\end{tabular}
\text{ with }a>b\text{ in }spl(T)\text{ and }%
\begin{tabular}
[c]{|l|l|l|l|}\hline
$.$ & $.$ & $b^{\prime}$ & $.$\\\hline
$.$ & $a^{\prime}$ & $.$ & $.$\\\hline
\end{tabular}
\text{ in }spl(T_{n}).\label{confs}%
\end{equation}

If $a^{\prime}\leq b^{\prime}$ we must have $a=a^{\prime}=b^{\prime}%
=\overline{n}$ and $b=n$ ($n$ is the last unbarred letter of $\mathcal{C}_{n}%
$) because the columns of $spl(T)$ and $spl(T_{n})$ differ by at most a change
$n\rightarrow\overline{n}$. We have $\overline{n}\in rC_{1}$, so $n\notin
C_{1}$ and $\overline{n}\in C_{1}$.\ Similarly, $n\notin D_{1}$ and
$\overline{n}\in D_{1}$. Hence $C_{1}=D_{1}$ for these two columns may only
differ by a change $n\rightarrow\overline{n}$.\ In this case, $\widetilde
{f}_{n}\,\mathrm{w}(C_{1}C_{2})=\mathrm{w(}C_{1})\,\widetilde{f}%
_{n}\mathrm{w(}C_{2})\neq0.$ Hence $n\in C_{2}$ and $\overline{n}\notin C_{2}%
$.\ By using the notation of (\ref{regle+-}), we have $\rho_{n}(\mathrm{w}%
(C_{1}C_{2}))=+-$, which implies that $\widetilde{f}_{n}\,\mathrm{w}(T)=0.$
Hence we obtain a contradiction.\ It means that $a^{\prime}>b^{\prime}$. So a
horizontal move occurs in $T_{n}$. Moreover, the letter $b^{\prime}$ slides
into $rD_{1}$ when we apply $\mathrm{SJDT}$ to $T_{n}$. If this were not so,
we would have the configurations
\[%
\begin{tabular}
[c]{|l|l|l|l|}\hline
$.$ & $.$ & $n$ & $.$\\\hline
$.$ & $n$ & $.$ & $.$\\\hline
\end{tabular}
\text{ in }T\text{ and }%
\begin{tabular}
[c]{|l|l|l|l|}\hline
$.$ & $.$ & $n$ & $.$\\\hline
$.$ & $\overline{n}$ & $.$ & $.$\\\hline
\end{tabular}
\text{ in }T_{n}%
\]
Hence $C_{1}$ would contain the letter $n$ which is changed in $\overline{n}$
when we apply $\widetilde{f}_{n}$ and $C_{2}\cap\{n,\overline{n}\}=\{n\}$.\ We
obtain a contradiction because we have $\rho_{n}(\mathrm{w}(C_{1}C_{2}))=--,$
hence by (\ref{regle+-}), it is the letter $n$ of $w(C_{2})$ that should be
changed in $\overline{n}$ when we apply $\widetilde{f}_{n}$.

When we compute $\mathrm{SJDT}$ on $T$ and $T_{n},$ we first execute a
sequence of vertical moves until we obtain a puncture $\ast$ $\ast$ on the row
containing $b$ and on the row containing $b^{\prime}$.\ Then we apply a
horizontal move on these two punctured tableaux.\ Denote by $C_{1}^{\prime
}C_{2}^{\prime}$ and $D_{1}^{\prime}D_{2}^{\prime}$ the skew punctured
tableaux so obtained.

Suppose first that $D_{1}=C_{1}$ and $\mathrm{w(}D_{2})=\widetilde{f}%
_{n}\,\mathrm{w}(C_{2})$. Then $D_{2}=C_{2}-\{n\}+\{\overline{n}\}$ and
$lD_{2}=lC_{2}-\{n\}+\{\overline{n}\}$ because $\overline{n}\notin C_{2}$. If
$b=b^{\prime}$ we can write $D_{1}^{\prime}D_{2}^{\prime}=C_{1}^{\prime}%
D_{2}^{\prime}$ and $D_{2}^{\prime}=C_{2}^{\prime}-\{n\}+\{\overline{n}\}$.
Then $\widetilde{f}_{n}\,\mathrm{w(}C_{1}^{\prime}C_{2}^{\prime}),$ which is
obtained from $\mathrm{w(}C_{1}^{\prime}C_{2}^{\prime})$ by turning $n$ into
$\overline{n}$ in $\mathrm{w(}C_{2}^{\prime}),$ is equal to $\mathrm{w(}%
C_{1}^{\prime}D_{2})$. When $b\neq b^{\prime},$ we have $b^{\prime}%
=\overline{n}$ and $b=n$. We can write $D_{2}^{\prime}=C_{2}^{\prime
}-\{n\}+\{\overline{n}\}$ and $D_{1}^{\prime}=C_{1}^{\prime}-\{n\}+\{\overline
{n}\}.$ Indeed $rC_{1},$ hence $C_{1},$ does not contain a letter of
$\{n,\overline{n}\}$. Then $\widetilde{f}_{n}\,\mathrm{w(}C_{1}^{\prime}%
C_{2}^{\prime}),$ which is obtained from $\mathrm{w(}C_{1}^{\prime}%
C_{2}^{\prime})$ by turning $n$ into $\overline{n}$ in $\mathrm{w(}%
C_{1}^{\prime}),$ is equal to $\mathrm{w(}D_{1}^{\prime}D_{2}^{\prime})$.

Now suppose that $D_{2}=C_{2}$ and $\widetilde{f}_{n}\,\mathrm{w}%
(C_{1})=\mathrm{w(}D_{1})$. We have $b=b^{\prime}$ hence $C_{2}^{\prime}%
=D_{2}^{\prime}$. Notice that $b\notin\{n,\overline{n}\}$ because
$\overline{n}\in rC_{1}$ and $n\in rD_{1}.$ So $D_{1}^{\prime}=C_{1}^{\prime
}-\{n\}+\{\overline{n}\}$. Then $\widetilde{f}_{n}\,\mathrm{w(}C_{1}^{\prime
}C_{2}^{\prime}),$ which is obtained from $\mathrm{w(}C_{1}^{\prime}%
C_{2}^{\prime})$ by turning $n$ into $\overline{n}$ in $\mathrm{w(}%
C_{1}^{\prime})$ is equal to $\mathrm{w(}D_{1}^{\prime}D_{2}^{\prime})$.
\end{proof}

\bigskip

The subalgebra of $U_{q}(sp_{2n})$ generated by the Chevalley generators
$e_{i},f_{i}$ and $t_{i}$ for $i=1,...,n-1$ is isomorphic to $U_{q}(sl_{n})$.
Each $U_{q}(sp_{2n})$-module $M$ is by restriction an $U_{q}(sl_{n}%
)$-module$.$ To obtain its crystal graph it suffices to erase the arrows of
color $n$ in the crystal graph of $M$. If $w\in G_{n}$, we denote by
$B^{A}(w)$ the $U_{q}(sl_{n})$-connected component of $G_{n}$ considered as a
$U_{q}(sl_{n})$-crystal graph (i.e. without the arrows of color $n$).

Set $\mathcal{U}_{(q,p)/k}^{+}=\{w\in\mathcal{U}_{(q,p)/k};$ $w$ contains only
unbarred letters$\}$ and $\mathcal{U}_{(q,p)/k}^{-}=\{w\in\mathcal{U}%
_{(q,p)/k};$ $w$ contains only barred letters$\}$. Then denote respectively by
$\Psi_{(q,p)/k}^{+}$ and $\Psi_{(q,p)/k}^{-}$ the restrictions of the map
$\Psi_{(q,p)/k}$ to $\mathcal{U}_{(q,p)/k}^{+}$ and $\mathcal{U}_{(q,p)/k}%
^{-}$.

\begin{lemma}
\label{lem_psi+}$\mathcal{U}_{(q,p)/k}^{+}$ and $\mathcal{U}_{(q,p)/k}^{-}$
are $U_{q}(sl_{n})$-sub-crystals of $G_{n}$ considered as a $U_{q}(sl_{n}%
)$-crystal.\ Moreover, $\Psi_{(q,p)/k}^{+}$ and $\Psi_{(q,p)/k}^{-}$ are
$U_{q}(sl_{n})$-crystal isomorphisms from $\mathcal{U}_{(q,p)/k}^{+}$ and
$\mathcal{U}_{(q,p)/k}^{-}$ to their images.
\end{lemma}

\begin{proof}
Barred and unbarred letters play a symmetric role for the action of
Kashiwara's operators $\widetilde{e}_{i},\widetilde{f}_{i}$ $i=1,...,n-1$. So
it suffices to prove the lemma for $\mathcal{U}_{(q,p)/k}^{+}$ and
$\Psi_{(q,p)/k}^{+}$.\ The set $\mathcal{U}_{(q,p)/k}^{+}$ is stable under the
action of $\widetilde{e}_{i},\widetilde{f}_{i}$ $i=1,...,n-1$ so it is a
$U_{q}(sl_{n})$-sub-crystal of $G_{n}$.\ Lemma \ref{lem_JT_HWV} show that
$\Psi_{(q,p)/k}^{+}$ induces a bijection from the set $\mathcal{U}%
_{(q,p)/k}^{0,+}$ of highest weight vertices of $\mathcal{U}_{(q,p)/k}^{0,+}$
to $\Psi_{(q,p)/k}^{+}(\mathcal{U}_{(q,p)/k}^{0,+}).$ Moreover if $w^{0}%
\in\mathcal{U}_{(q,p)/k}^{0,+}$, $\Psi_{(q,p)/k}(w^{0})$ and $w^{0}$ have the
same weight. Then $\Psi_{(q,p)/k}^{+}$ is a $U_{q}(sl_{n})$-crystal
isomorphism if and only if for any $w\in\mathcal{U}_{(q,p)/k}^{+}$%
\begin{equation}
\widetilde{f}_{i}\circ\Psi_{(q,p)/k}^{+}(w)=\Psi_{(q,p)/k}^{+}\circ
\widetilde{f}_{i}(w)\text{ }i=1,...,n-1.\label{JT_comp_ei}%
\end{equation}
This equality was proved in Theorem 3.3.1 of \cite{VL}.
\end{proof}

\bigskip

Consider $r$ and $s$ two integers such that $r+s\leq n$.\ The vertices
$1\cdot\cdot\cdot r\overline{n}\cdot\cdot\cdot(\overline{n-s+1})$ and
$\overline{n}\cdot\cdot\cdot(\overline{n-s+1})1\cdot\cdot\cdot r$ are highest
weight vertices with the same $U_{q}(sl_{n})$-weight of the crystal $G_{n}$
considered as a $U_{q}(sl_{n})$-crystal graph.\ Denote by $\Theta_{r,s}$ the
$U_{q}(sl_{n})$-crystal isomorphism:%
\[
B^{A}(1\cdot\cdot\cdot r\overline{n}\cdot\cdot\cdot(\overline{n-s+1}%
))\overset{U_{q}(sl_{n})}{\rightarrow}B^{A}(\overline{n}\cdot\cdot
\cdot(\overline{n-s+1})1\cdot\cdot\cdot r).
\]
Note that the vertices of $B^{A}(1\cdot\cdot\cdot r\overline{n}\cdot\cdot
\cdot(\overline{n-s+1}))$ are admissible column words of length $r+s$ with $r$
unbarred letters and those of $B^{A}(\overline{n}\cdot\cdot\cdot
(\overline{n-s+1})1\cdot\cdot\cdot r)$ are words of the form $u^{-}v^{+}$
where the factors $u^{-}$ and $v^{+}$ contain respectively $s$ barred and $r$
unbarred letters. Then we have the following

\begin{lemma}
\label{lem_theta_comm_with_ei} \ \ \ \ \ \ \ \ \ \ \ \ \ 

\begin{enumerate}
\item $B^{A}(1\cdot\cdot\cdot r\overline{n}\cdot\cdot\cdot(\overline{n-s+1}))$
consists of the readings of the admissible column words of length $r+s$ with
$r$ unbarred letters.

\item $B^{A}(\overline{n}\cdot\cdot\cdot(\overline{n-s+1})1\cdot\cdot\cdot r)$
consists of the words of the form $u^{-}v^{+}$ such that $v^{+}u^{-}$ is the
reading of a coadmissible column of height $r+s$ with $r$ unbarred letters.

\item  Let $C$ be a column such that $\mathrm{w(}C)\in B^{A}(1\cdot\cdot\cdot
r\overline{n}\cdot\cdot\cdot(\overline{n-s+1}))$.\ Write $\mathrm{w(}C^{\ast
})=v^{+}u^{-}$ where the words $u^{-}$ and $v^{+}$ contain respectively barred
and unbarred letters.\ Then%
\[
\Theta_{r,s}(\mathrm{w(}C))=u^{-}v^{+}.
\]
\end{enumerate}
\end{lemma}

\begin{proof}
$1.\;$Follows from the fact that $w_{r,s}=1\cdot\cdot\cdot r\overline{n}%
\cdot\cdot\cdot(\overline{n-s+1})$ is the unique admissible column word of
length $r+s$ with $r$ unbarred letters satisfying $\varepsilon_{i}(w_{r,s})=0$
for $i=1,...,n-1$. In the sequel we set $w_{r,s}^{\ast}=\overline{n}\cdot
\cdot\cdot(\overline{n-s+1})1\cdot\cdot\cdot r$.

$2.$ Consider $w\in B^{A}(w_{r,s}^{\ast})$ and suppose that for any pair
$(z,\overline{z})\in w$
\begin{equation}
N^{\ast}(z)\text{ the number of letters }x\in w\text{ such }z\leq
x\leq\overline{z}\text{ verifies }N^{\ast}(z)\leq n-z+1.\label{cond_N_star}%
\end{equation}
Let $i\in\{1,...,n-1\}$ such that $\widetilde{f}_{i}(w)\neq0$.\ Then any pair
of letters $(z,\overline{z})\in\widetilde{f}_{i}(w)$ verifies $N^{\ast}(z)\leq
n-z+1$.\ Otherwise $\widetilde{f}_{i}(w)$ contains a pair $(t,\overline{t})$
with $N^{\ast}(t)>n-t+1$.\ Then $t\notin w$ or $\overline{t}\notin w$.\ By
(\ref{regle+-}) we have $i=t$ and $w$ is obtained from $\widetilde{f}_{i}(w)$
by turning $\overline{t}$ into $\overline{t+1}$.\ Then the number of letters
$x\in w$ such that $t\leq x\leq\overline{t+1}$ is $>n-t+1$.\ Let $y$ be the
smallest letter of $w$ such that $t<y\leq n$ and $(y,\overline{y})\in w$.\ We
have $N^{\ast}(y)>n-y+1$.\ Hence we derive a contradiction. It is clear that
$w_{r,s}^{\ast}$ verifies (\ref{cond_N_star}).\ Hence by induction all the
vertices of $B^{A}(w_{r,s}^{\ast})$ verify (\ref{cond_N_star}).\ By
(\ref{N_star}), it means that these vertices may be written $u^{-}v^{+}$ where
$v^{+}u^{-}$ is the reading of a coadmissible column of height $r+s$ with $r$
unbarred letters. It follows from (\ref{phi}) that the number of admissible
columns of height $r+s$ with $r$ unbarred letters is equal to the number of
coadmissible columns of height $r+s$ with $r$ unbarred letters. So $2$ is proved.

$3.$ For any $w=\mathrm{w(}C)\in B^{A}(w_{r,s})$, we set $w^{\ast}=u^{-}v^{+}$
where $\mathrm{w}(C^{\ast})=v^{+}u^{-}$.\ By $1$, $2$ and (\ref{phi}) we know
that the map $w\rightarrow w^{\ast}$ is a bijection from $B^{A}(w_{r,s})$ to
$B^{A}(w_{r,s}^{\ast})$.\ Note that $\Theta_{r,s}(w_{r,s})=w_{r,s}^{\ast}%
.\;$Hence, to prove that $\Theta_{r,s}(w)=w^{\ast}$, it suffices to show that
for any $i=1,...,n-1$ such that $\widetilde{f}_{i}(w)\neq0:$%
\begin{equation}
(\widetilde{f}_{i}(w))^{\ast}=\widetilde{f}_{i}(w^{\ast}).\label{star_com_fi}%
\end{equation}
Fix $i\in\{1,...,n-1\}$ and set $E_{i}=w\cap\{i,i+1,\overline{i+1}%
,\overline{i}\},$ $F_{i}=w^{\ast}\cap\{i,i+1,\overline{i+1},\overline{i}%
\}.\;$Then by (\ref{regle+-}) $E_{i}$ is equal to one of the following sets:
$\mathrm{(i)}$ $E_{i}=\{i\}$, $\mathrm{(ii)}$ $E_{i}=\{\overline{i+1}\},$
$\mathrm{(iii)}$ $E_{i}=\{i,\overline{i+1},\overline{i}\},$ $\mathrm{(iv)}$
$E_{i}=\{i,i+1,\overline{i+1}\}$, $\mathrm{(v)}$ $E_{i}=\{i+1,\overline
{i+1}\}$ or $\mathrm{(vi)}$ $E_{i}=\{i,\overline{i+1}\}.$\ Then
(\ref{star_com_fi}) follows by considering for the six cases above, the
possible sets $F_{i}$.\ For example in case $\mathrm{(i)}$, we have
$\widetilde{f}_{i}(w)=w-\{i\}+\{i+1\}$. Moreover $F_{i}=\{i\}$ or
$F_{i}=\{i,i+1,\overline{i+1}\}$.\ Then $(\widetilde{f}_{i}(w))^{\ast
}=\widetilde{f}_{i}(w^{\ast})=w^{\ast}-\{i\}+\{i+1\}$ when $F_{i}=\{i\}$ and
by (\ref{regle+-}), $(\widetilde{f}_{i}(w))^{\ast}=\widetilde{f}_{i}(w^{\ast
})=w^{\ast}-\{\overline{i+1}\}+\{\overline{i}\}$ when $F_{i}=\{i,i+1,\overline
{i+1}\}$.
\end{proof}

\bigskip

The proof of Proposition \ref{Prop_phsi_iso} relies on the following simple

\begin{lemma}
\label{lem_util}Let $\Gamma_{i}$ and $\Gamma_{i}^{\prime}$ $i=1,...,m$ be
sub-crystals of $G_{n}$ such that there exists a crystal isomorphism $\xi_{i}$
from $\Gamma_{i}$ to $\Gamma_{i}^{\prime}$ for $i=1,...,m$. Denote
respectively by $\Gamma_{1}\cdot\cdot\cdot\Gamma_{m}$ and $\Gamma_{1}^{\prime
}\cdot\cdot\cdot\Gamma_{m}^{\prime}$ the sub-crystals of $G_{n}$ whose
vertices are the words of the form $w_{1}\cdot\cdot\cdot w_{m}$ obtained by
concatenating the words $w_{i}\in\Gamma_{i}$ $i=1,...,m$ \textrm{(}%
resp.\ $w_{i}\in\Gamma_{i}^{\prime}$ $i=1,...,m$\textrm{)}.\ Then the map%
\begin{gather*}
\Gamma_{1}\cdot\cdot\cdot\Gamma_{m}\rightarrow\Gamma_{1}^{\prime}\cdot
\cdot\cdot\Gamma_{m}^{\prime}\\
w_{1}\cdot\cdot\cdot w_{m}\longmapsto\xi_{1}(w_{1})\cdot\cdot\cdot\xi
_{m}(w_{m})
\end{gather*}
is a crystal isomorphism.
\end{lemma}

\begin{proof}
This follows immediately from the description (\ref{tens}) of the tensor
product of two crystal graphs.
\end{proof}

Now we can state the

\begin{proposition}
\label{Prop_phsi_iso}The map $\Psi_{(q,p)/k}$ is a crystal isomorphism from
$\mathcal{U}_{(q,p)/k}$ to its image.
\end{proposition}

\begin{proof}
By Lemmas \ref{lem_iso_G}, \ref{lem_JT_HWV} and \ref{lem_com_en}, it suffices
to prove that:%
\begin{equation}
\widetilde{f}_{i}\circ\Psi_{(q,p)/k}(w)=\Psi_{(q,p)/k}\circ\widetilde{f}%
_{i}(w)\text{ for }i=1,...,n-1.\label{phsi_com_fi}%
\end{equation}
This will be deduced from Lemmas \ref{lem_psi+}, \ref{lem_theta_comm_with_ei}
and \ref{lem_util}.\ Set $w=\mathrm{w}(C_{1}C_{2})$ with $C_{1}C_{2}%
\in\mathcal{T}_{(q,p)/k}$. If no horizontal move occurs in $C_{1}C_{2},$
(\ref{phsi_com_fi}) is immediate.\ Otherwise, suppose that an unbarred letter
of $lC_{2}$ slides into $rC_{1}$ during the horizontal move.\ Write
$\mathrm{w(}C_{1})=u_{1}^{+}u_{1}^{-}$ and $\mathrm{w(}C_{2})=u_{2}^{+}%
u_{2}^{-}$ where the letters of $u_{1}^{+}$ and $u_{2}^{+}$ are unbarred and
those of $u_{1}^{-}$ and $u_{2}^{-}$ are barred. Denote respectively by
$r_{1},r_{2}$ the number of unbarred letters of $u_{1}^{+}$ and $u_{2}^{+}%
$.\ Then $u_{1}^{-}$ and $u_{2}^{-}$ contain respectively $q-k-r_{1}=s_{1}$
and $p-r_{2}=s_{2}$ barred letters.\ Set $\Theta_{r_{2},s_{2}}(u_{2}^{+}%
u_{2}^{-})=v_{2}^{-}v_{2}^{+}$ (see Lemma \ref{lem_theta_comm_with_ei}). We
can decompose the computation of $\Psi_{(q,p)/k}\,(\mathrm{w(}C_{1}%
C_{2}))=\Psi_{(q,p)/k}(\,u_{2}^{+}u_{2}^{-}$ $u_{1}^{+}u_{1}^{-})$ into the
following three steps (as illustrated by Example \ref{ex_illustration}):

\noindent$\mathrm{(i)}$ we calculate the word%
\[
v_{2}^{-}v_{2}^{+}u_{1}^{+}u_{1}^{-},
\]

\noindent$\mathrm{(ii)}$ we apply $\Psi_{(k+r_{1},r_{2})/k}^{+}$ to the
admissible skew tableau $v_{2}^{+}u_{1}^{+}$ containing only unbarred letters
to obtain%
\[
v_{2}^{-}\,\Psi_{(k+r_{1},r_{2})/k}^{+}(v_{2}^{+}u_{1}^{+})\,u_{1}^{-}%
=v_{2}^{-}\,b_{2}^{+}b_{1}^{+}\,u_{1}^{-}%
\]
where $b_{1}^{+}$ and $b_{2}^{+}$ contain respectively $r_{1}+1$ and $r_{2}-1$
unbarred letters.

\noindent$\mathrm{(iii)}$ we apply $\Theta_{r_{2}-1,s_{2}}^{-1}$ to $v_{2}%
^{-}b_{2}^{+},$ getting%
\[
\Psi_{(q,p)/k}\,\mathrm{w(}C_{1}C_{2})=\Theta_{r_{2}-1,s_{2}}^{-1}(v_{2}%
^{-}b_{2}^{+})\,b_{1}^{+}u_{1}^{-}.
\]

By Lemmas \ref{lem_psi+} and \ref{lem_theta_comm_with_ei}, the maps
$\Psi_{(k+r_{1},r_{2})/k}^{+},$ $\Theta_{r_{2}-1,s_{2}}^{-1}$ and
$\Theta_{r_{2},s_{2}}$ are $U_{q}(sl_{n})$-crystal isomorphisms. Hence
$\Psi_{(q,p)/k}$ commutes with the operators $\widetilde{f}_{i}$
$i=1,...,n-1$. This follows from Lemma \ref{lem_util} applied with $m=2,$
$\xi_{1}=\Theta_{r_{2},s_{2}}$, $\xi_{2}=id$ in step $\mathrm{(i)}$; with
$m=3,$ $\xi_{1}=\xi_{3}=id$, $\xi_{2}=\Psi_{(k+r_{1},r_{2})/k}^{+}$ in step
$\mathrm{(ii)}$ and with $m=2,$ $\xi_{1}=\Theta_{r_{2}-1,s_{2}}^{-1}$,
$\xi_{2}=id$ in step $\mathrm{(iii)}$.

When a barred letter of $lC_{2}$ slides into $rC_{1},$ we obtain that
$\Psi_{(q,p)/k}$ commute with the operators $\widetilde{f}_{i}$ $i=1,...,n-1$
by similar arguments. This time we compute $\Theta_{r_{1},s_{1}}(u_{1}%
^{+}u_{1}^{-})$ and we use the restriction $\Psi_{(e+s_{1},s_{2})/e}^{-}$ with
$e=k+r_{1}-r_{2}$ instead of $\Psi_{(k+r_{1},r_{2})/k}^{+}$.
\end{proof}

\begin{example}
\label{ex_illustration}Consider $C_{1}C_{2}$ such that $spl(C_{1}C_{2})=%
\begin{tabular}
[c]{|l|l|ll}\hline
$\mathtt{\ast}$ & $\mathtt{\ast}$ & $\mathtt{2}$ &
\multicolumn{1}{|l|}{$\mathtt{2}$}\\\hline
$\mathtt{2}$ & $\mathtt{2}$ & $\mathtt{3}$ & \multicolumn{1}{|l|}{$\mathtt{3}%
$}\\\hline
$\mathtt{3}$ & $\mathtt{3}$ & $\mathtt{4}$ & \multicolumn{1}{|l|}{$\mathtt{5}%
$}\\\hline
$\mathtt{\bar{5}}$ & $\mathtt{\bar{5}}$ & $\mathtt{\bar{5}}$ &
\multicolumn{1}{|l|}{$\mathtt{\bar{4}}$}\\\hline
$\mathtt{\bar{1}}$ & $\mathtt{\bar{1}}$ &  & \\\cline{1-2}%
\end{tabular}
$.\ We have $\mathrm{w(}C_{1})=23\bar{5}\bar{1}$ and $\mathrm{w(}%
C_{2})=235\bar{5}$ so $\mathrm{w(}C_{1}C_{2})=235\bar{5}\,23\bar{5}\bar{1}%
$.\ When we apply $\mathrm{SJDT}$ to $C_{1}C_{2},$ the letter $4\in lC_{2}$
slides into $rC_{1}$ (see the punctured skew tableau $T_{2}$ of Example
\ref{ex_horizontal_move}).\ With the notation of the above proof we have
$v_{2}^{-}v_{2}^{+}u_{1}^{+}u_{1}^{-}=\bar{4}234\,23\bar{5}\bar{1}$.\ By
considering the slide
\[%
\begin{tabular}
[c]{|l|l|}\hline
$\mathtt{2}$ & $\mathtt{2}$\\\hline
$\mathtt{3}$ & $\mathtt{3}$\\\hline
$\mathtt{\ast}$ & $\mathtt{4}$\\\hline
\end{tabular}
\rightarrow%
\begin{tabular}
[c]{|l|l|}\hline
$\mathtt{2}$ & $\mathtt{2}$\\\hline
$\mathtt{3}$ & $\mathtt{3}$\\\hline
$\mathtt{4}$ & $\mathtt{\ast}$\\\hline
\end{tabular}
\]
we obtain $\Psi_{(3,3)/1}^{+}(234\,23)=23\,234$.\ Then $v_{2}^{-}\,b_{2}%
^{+}b_{1}^{+}\,u_{1}^{-}=\bar{4}23\,234\bar{5}\bar{1}$ and $\Psi
_{(5,4)/1}\,\mathrm{w(}C_{1}C_{2})=23\bar{4}\,2345\bar{1}$.\ This last word is
the reading of
\begin{tabular}
[c]{|l|l}\hline
$\mathtt{2}$ & \multicolumn{1}{|l|}{$\mathtt{2}$}\\\hline
$\mathtt{3}$ & \multicolumn{1}{|l|}{$\mathtt{3}$}\\\hline
$\mathtt{4}$ & \multicolumn{1}{|l|}{$\mathtt{\bar{4}}$}\\\hline
$\mathtt{\bar{5}}$ & \\\cline{1-1}%
$\mathtt{\bar{1}}$ & \\\cline{1-1}%
\end{tabular}
the skew tableau obtained by applying $\mathrm{SJDT}$ on $C_{1}C_{2}.$
\end{example}

Note that the sub-crystal $\Psi_{(q,p)/k}(\mathcal{U}_{(q,p)/k})$ may contain
readings of skew tableaux with different shapes. We can now state the

\begin{theorem}
\label{th_iso_JDT}Let $c$ be a fixed inner corner of the shape $(\lambda/\mu
)$.\ Denote by $\Xi(\cdot,c)$ the map defined on $\mathcal{U}_{(\lambda/\mu)}$
by $\Xi(w,c)=\mathrm{w(SJDT(}T_{w},c))$ where $T_{w}$ is the admissible skew
tableau of $\mathcal{T}_{(\lambda/\mu)}$\ of reading $w$. Then $\Xi(\cdot,c)$
is a crystal isomorphism from $\mathcal{U}_{(\lambda/\mu)}$ to its image.
\end{theorem}

\begin{proof}
By Proposition \ref{Prop_phsi_iso} and Lemma \ref{lem_util}, each step of
$\mathrm{SJDT(}T_{w},c)$ may be interpreted on the readings of intermediate
tableaux as the result of the action of a crystal isomorphism. Hence
$\Xi(\cdot,c)$ is a crystal isomorphism.
\end{proof}

\begin{corollary}
Let $T$ be an admissible skew tableau. Then by applying the $\mathrm{SJDT}$
successively to the inner corners of $T$ we obtain a symplectic tableau
independent of the order in which these inner corners are filled. Moreover
this tableau coincides with $P(\mathrm{w(}T))$.
\end{corollary}

\begin{proof}
By the previous theorem, all the possible readings $w$ of symplectic tableaux
obtained from $T$ by iterating the $\mathrm{SJDT}$ must satisfy $w\sim
\mathrm{w}(T)$. Hence $w=\mathrm{w}(P(T))$ is independent of the order in
which the inner corners are filled.
\end{proof}

\bigskip

\noindent\textbf{Acknowledgments:} The author thanks the anonymous referee for
his very careful reading of the paper and his many valuable corrections and comments.

\end{document}